
\parindent=0pt\magnification=1200\font\third=cmbx10 at 15pt

\font\arrfont=cmbsy10
\newfam\nwarrowfamm

\textfont\nwarrowfamm=\arrfont
\def\nwarrow{{\arrfont\char45}}
\def\li{\l}

\def\Pn#1{{\bf P}^#1}
\def\Pnd#1{{\check{\bf P}}^#1}
\def\es#1#2#3#4#5#6{0\longrightarrow \; {\cal O}_{#1}(#2)\; \longrightarrow \;
{\cal O}_{#3}(#4)\;
\longrightarrow \; {\cal O}_{#5}(#6) \; \longrightarrow 0}

\def\esi#1#2#3#4#5#6{0\longrightarrow \;  {\cal I}_{#1}(#2)\; \longrightarrow
\;
{\cal I}_{#3}(#4)\; \longrightarrow \; {\cal I}_{#5}(#6) \; \longrightarrow 0}
\def\eso#1#2#3#4#5#6{0\longrightarrow \;  {\cal I}_{#1}(#2)\; \longrightarrow
\;
{\cal O}_{#3}(#4)\; \longrightarrow \; {\cal O}_{#5}(#6) \; \longrightarrow 0}
\def\id#1#2{{\cal I}_{#1}(#2)}
\def\kn#1#2{{\cal O}_{#1}(#2)}
\def\phi#1{\varphi_{#1}}
\def\coh#1#2#3{h^{#1} (\kn{#2}{#3})}
\def\Coh#1#2#3{{\rm H}^{#1} ({\kn{#2}{#3}})}
\def\cohi#1#2#3{h^{#1} ({\id{#2}{#3}})}

\def\slutt{\sqcap\!\!\!\!\sqcup}
\def\rto{\raise.5ex\hbox{$\scriptscriptstyle ---\!\!\!>$}}

\chardef\po="19
\def\deh#1{{\Delta}({#1})}


\def\li{\l}

\def\Pn#1{{\bf P}^#1}
\def\Pnd#1{{\check{\bf P}}^#1}
\def\es#1#2#3#4#5#6{0\longrightarrow \; {\cal O}_{#1}(#2)\; \longrightarrow \;
{\cal O}_{#3}(#4)\;
\longrightarrow \; {\cal O}_{#5}(#6) \; \longrightarrow 0}

\def\esi#1#2#3#4#5#6{0\longrightarrow \;  {\cal I}_{#1}(#2)\; \longrightarrow
\;
{\cal I}_{#3}(#4)\; \longrightarrow \; {\cal I}_{#5}(#6) \; \longrightarrow 0}
\def\eso#1#2#3#4#5#6{0\longrightarrow \;  {\cal I}_{#1}(#2)\; \longrightarrow
\;
{\cal O}_{#3}(#4)\; \longrightarrow \; {\cal O}_{#5}(#6) \; \longrightarrow 0}
\def\ess#1#2#3#4#5#6#7#8{0\longrightarrow \; {\cal O}_{#1}(#2)\;
\longrightarrow \; {\cal O}_{#3}(#4)\;
\longrightarrow \; {\cal O}_{#5}(#6)\oplus {\cal O}_{#7}(#8) \; \longrightarrow
0}

\def\id#1#2{{\cal I}_{#1}(#2)}
\def\kn#1#2{{\cal O}_{#1}(#2)}
\def\phi#1{\varphi_{#1}}
\def\coh#1#2#3{h^{#1}(\kn {#2}{#3})}
\def\Coh#1#2#3{{\rm H}^{#1}(\kn{#2}{#3})}
\def\cohi#1#2#3{h^{#1}(\id {#2}{#3})}

\def\slutt{\sqcap\!\!\!\!\sqcup}
\def\rto{\raise.5ex\hbox{$\scriptscriptstyle ---\!\!\!>$}}

\chardef\po="19
\def\deh#1{{\Delta}({#1})}

\def\mapright#1{\smash{\mathop{\longrightarrow}\limits^{#1}}}
\def\dsum{\mathop{\oplus}}
\def\maponto#1{\smash{\mathop{\longrightarrow}\limits^{#1}}
\hskip -0.2cm{\raise.3ex\hbox{$\scriptscriptstyle >$}}}
\def\mapontoa#1{\smash{\mathop{\longrightarrow}\limits^{#1}}
\hskip -0.07cm{\raise.3ex\hbox{$\scriptscriptstyle >$}}}

\def\struct#1#2{{\cal O}_{#1}(#2)}
\def\sstruct#1{{\cal O}_{#1}}
\def\openP{\bf P}
\def\hi#1{{h}^{#1}}
\def\scohid#1#2#3{{ h}^{#1}({\cal I}_{#2}({#3}))}
\def\scohos#1#2#3{{ h}^{#1}({\cal O}_{#2}({#3}))}
\def\sCohid#1#2#3{{ H}^{#1}({\cal I}_{#2}({#3}))}

\def\exactg#1#2#3{0\longrightarrow \;{#1}\; \longrightarrow \;
{#2}\; \longrightarrow \;  {#3}\; \longrightarrow 0}


\newbox\ersttest
\newdimen\breite
\breite=0pt
\def\test#1#2#3#4#5#6#7#8{
\setbox\ersttest=\hbox{#1}
\ifnum \breite<\wd\ersttest
  \breite=\wd\ersttest
\else
\fi
\setbox\ersttest=\hbox{#2}
\ifnum \breite<\wd\ersttest
  \breite=\wd\ersttest
\else
\fi
\setbox\ersttest=\hbox{#3}
\ifnum \breite<\wd\ersttest
  \breite=\wd\ersttest
\else
\fi
\setbox\ersttest=\hbox{#4}
\ifnum \breite<\wd\ersttest
  \breite=\wd\ersttest
\else
\fi
\setbox\ersttest=\hbox{#5}
\ifnum \breite<\wd\ersttest
  \breite=\wd\ersttest
\else
\fi
\setbox\ersttest=\hbox{#6}
\ifnum \breite<\wd\ersttest
  \breite=\wd\ersttest
\else
\fi
\setbox\ersttest=\hbox{#7}
\ifnum \breite<\wd\ersttest
  \breite=\wd\ersttest
\else
\fi
\setbox\ersttest=\hbox{#8}
\ifnum \breite<\wd\ersttest
  \breite=\wd\ersttest
\else
\fi
}
\long\def\beili55v#1#2#3#4#5#6#7#8#9{
\vbox{\tabskip=0pt\offinterlineskip
\def\tablerule{\omit&\multispan{11}{\hrulefill}&&\cr}
\halign to250pt{\strut##  & \vrule##\tabskip=1em  &\hbox
to\breite{\hfil##\hfil}&\vrule##&
              \hbox to\breite{\hfil##\hfil}& \vrule##&\hbox
to\breite{\hfil##\hfil}&\vrule##&
              \hbox to\breite{\hfil##\hfil}& \vrule##&\hbox
to\breite{\hfil##\hfil}&\vrule##&
              \hbox to\breite{\hfil##\hfil}& ##\hfil &\hfil## \tabskip=0em \cr
\hskip-13pt
i&\omit\hskip-2.2pt$\uparrow$&&\omit&&\omit&&\omit&&\omit&&\omit&&&\cr
\noalign{\vskip-2pt}
&\omit\vrule height10pt&\omit &\omit&&\omit&&\omit&&\omit&&\omit&&&\cr%
\tablerule
\omit&height5pt&\omit&&\omit&&\omit&&\omit&&\omit&&\omit&\omit&\cr
&& && && && && &&&&\cr
\tablerule
\omit&height5pt&\omit&&\omit&&\omit&&\omit&&\omit&&\omit&\omit&\cr
&&#1&& && && && &&&&\cr
\tablerule
\omit&height5pt&\omit&&\omit&&\omit&&\omit&&\omit&&\omit&\omit&\cr
&&#2&&#3&&#4&& && &&&\quad\quad\quad#9&\cr
\tablerule
\omit&height5pt&\omit&&\omit&&\omit&&\omit&&\omit&&\omit&\omit&\cr
&& && &&#5&&#6&&#7&&&&\cr
\tablerule
\omit&height5pt&\omit&&\omit&&\omit&&\omit&&\omit&&\omit&\omit&\cr
&& && && && &&#8&&&&\cr
\noalign{\vskip-3.1pt}
\omit&\multispan{12}{\rightarrowfill}&\cr
\noalign{\vskip5pt}
\multispan{13}{\hfil p}&\cr
}
}
}

{\hbox {}}\bigskip
\centerline{\third   Surfaces of Degree 10
in the Projective Fourspace}

\bigskip
\centerline{\third  via Linear systems and Linkage}

\bigskip \bigskip \centerline{\bf Sorin Popescu and Kristian Ranestad }
\bigskip\bigskip\bigskip
{\bf Contents}\bigskip
{\settabs 3\columns \+0. Introduction     &&1\cr\+1. 6-secants and plane
curves    &&6\cr\+2.  Linear systems     &&10\cr\+3. Syzygies    &&31\cr\+4.
Linkage     &&39\cr}
\bigskip\bigskip

{\bf 0\quad Introduction} \bigskip
Ellingsrud and Peskine showed in [ElP] that there are finitely many components
of
the Hilbert scheme of $\Pn 4$ containing smooth surfaces not of general type.
The upper bound for the degree of such surfaces has recently been reduced to
105
(cf. [BF]). Quite a bit of work has been put into trying to construct such
surfaces
of high degree.   So far the record is 15, which coincides with the conjectural
upper bound.
This paper concerns the classification of surfaces of degree 10 and
sectional genus 9 and 10. The surfaces of degree at most 9 are described
through classical work dating from the last century up to recent years:
[Ba], [Ro], [Io], [Ok], [Al], [AR]. Surfaces of degree
at least 11 have been considered systematically recently in [Po]. In degree
10 there are the abelian surfaces discovered by Comessatti [Co], rediscovered
by Horrocks and
Mumford [HM] in the seventies as zero-sections of an indecomposable rank 2
vector bundle on
$\Pn 4$. Beside surfaces of sectional genus at least 11, which can be linked to
smooth surfaces of smaller degree, other surfaces of degree 10 were considered
only recently.
Serrano gave examples of bielliptic surfaces of degree 10 in [S], which have
sectional genus 6 like
the abelian ones.  The Hilbert schemes of abelian and bielliptic surfaces are
 now well understood (cf. [BHM], [BaM], [H], [HKW], [HL], [HV], [L], [R],
[ADHPR]).
The second author determined numerical invariants and gave some examples of
surfaces with
sectional genus 8,9 and 10.  There are two families of surfaces of genus 8, one
of
rational surfaces, and one of non-minimal
Enriques surfaces. A construction using syzygies of both of them is given in
[DES].  The linear system of surfaces of the first family is described in [Ra],
cf. also [Al],
while for the second it is described in [Br]. The purpose of this paper is to
describe the remaining components
of the Hilbert scheme, i.e., to describe the smooth surfaces in $\Pn 4$ of
degree 10 and sectional genus 9 and 10. The results in themselves are to be
considered
interesting from the perspective of the diversity of techniques with which we
present them.
Thus we use relations between multisecants, linear systems, syzygies and
linkage to
describe the geometry of each surface.  We want in fact to stress the
importance of multisecants and syzygies for the study of these surfaces.
Adjunction which
provided efficient arguments for the classification of surfaces of smaller
degrees, here appears
to be less effective and will play almost no role in the proofs.
\medskip

  We show that there are 8 different families of smooth surfaces of degree
10 and sectional genus 9 and 10. The families are determined by numerical data
such as the
sectional genus $\pi$, the Euler characteristic $\chi =\chi ({\cal O}_S)$,
the number $N_6$ of 6-secants to the surface $S$ and the number $N_5$ of
5-secants to the surface which meet a general plane.  For each type we
describe the linear system of hyperplane sections on $S$, the resolution of
the ideal, the geometry of the surface in terms of curves on the surface
and hypersurfaces containing the surface, and the liaison class; in particular
minimal elements in the even liaison class.  Each type corresponds to an
irreducible component of the Hilbert scheme, and the dimension is computed.
The following table collects the numerical data and give references to the
points in this paper where the corresponding additional information can be
found.

$$\matrix {S&\pi &\chi &N_6 &N_5 & {\rm Cohomology}& {\rm Birational} &
{\rm Linear}  &{\rm construc-} &{\rm linkage} & {\rm Hilbert}\cr &&&&&&
{\rm type}&{\rm system}&{\rm tion}&&{\rm
scheme }\cr\cr A&9&1&1&6&(3.4)&{\rm rational}&(2.1)&(4.5)&(4.3)& 42\cr
B&9&1&0&6&(3.3)&{\rm rational}&(2.1)&(2.10)&(4.2)&42\cr
C&9&2&3&12&(3.7)&{\rm K3}&(2.12)&(2.14)&(4.7)&45\cr
D&9&2&1&12&(3.6)&{\rm K3}&(2.19)&(4.10)&(4.9)&44\cr
E&9&2&0&12&(3.5)&{\rm elliptic}&(2.23)&(2.31)&(4.12)&44\cr
F&9&3&3&18&(3.9)&{\rm gen. type}&(2.32)&(2.33)&(4.13)&47\cr
G&10&3&0&2&(3.11)&{\rm elliptic}&(2.34)&(4.15)&(4.14)&51\cr
H&10&4&1&6&(3.12)&{\rm gen. type}&(2.34)&(4.20)&(4.18)&53\cr}$$
\smallskip

The organization of the paper follows the different approaches to a
description
of the surfaces rather than each surface one by one.  Thus we focus on the
different methods used.  In the first section we use 6-secants and plane
curves
to get information on the Hartshorne-Rao module of the surface.  In the
second
section we use a geometric approach to describe linear systems and special
curves on $S$ and special hypersurfaces containing $S$.  In the third
section
the constructions via the Eagon-Northcott complex performed in [DES] is
recovered together with the resolution of the ideal of $S$. In the fourth
section we use both the resolution of the ideal and the geometry of the
surface to describe the minimal element in the even liaison class of
$S$.\bigskip

   {\bf Notation and basic results:}\par
\proclaim Adjunction formula 0.1.
$2p_a(C)-2=C^2+C\cdot K$.\par $Proof.$ See [Ha, Prop. 1.5].$\slutt$\bigskip
For curves $C$, $D$ and $C\cup D$ on a smooth surface $S$ the adjunction
immediately gives the following addition formula for the arithmetic genus:
$$p_a(C\cup D)=
p_a(C)+p_a(D)+C\cdot D-1.\leqno (0.2.)$$ This formula quickly yields the
following

\proclaim Lemma 0.3.  Let $C$ be a non-planar curve of
degree $d$ and arithmetic genus $p$ on a smooth surface. If $d=5$, then $p\leq
3$ with equality
only if $C$ decomposes into a plane quartic and a line meeting in a point.
When $d=6$, then $5\leq p\leq 6$  only if $C$ decomposes into a plane
quintic curve and a line, while $3\leq p\leq 4$ implies that $C$ spans a $\Pn
3$
unless $p=3$ and $C$ decomposes into a plane quartic curve and a conic, or two
lines, meeting the quartic in a point.\par
$Proof.$  Straightforward from (0.2) and the genus bound for
irreducible curves.$\slutt$

 \proclaim  Theorem (Riemann-Roch) 0.4.
$$\chi (\kn SC)=\coh 0SC-\coh 1SC+\coh 0S{K-C}={1\over 2} (C^2-C\cdot K)
+\chi (S).$$\par $Proof.$ See [Ha, Th.1.6].$\slutt$

\proclaim Hodge index theorem 0.5. If $H$ is an ample divisor and $D$ is a
divisor on $S$ such
that $H\cdot D=0$, then either $D^2<0$, or $D$ is numerically equivalent to
$0$.\par
$Proof.$ See [Ha, Th.1.9].$\slutt$\par\bigskip

For smooth surfaces in $\Pn 4$ with normal
bundle $N_S$ there is the relation, $$d^2-c_2(N_S)=d^2-10d-5H\cdot
K-2K^2+12\chi(S)=0,\leqno(0.6.)$$ which expresses the fact that $S$ has no
double points.  This will be referred to as the double point formula.

\proclaim Theorem (Severi) 0.7. All smooth surfaces in $\Pn 4$, except for the
Veronese surfaces, are linearly normal.\par $Proof.$ [Se] or
[Mo].$\slutt$\bigskip

Some classical numerical formulae for multisecant lines to a smooth surface in
$\Pn 4$ have recently been studied again by Le Barz: \bigskip

{\bf Multisecants 0.8.} ([LB]) If $\sharp_5$ denotes the number of 5-secant
lines to $S$ which meets a general plane, and $\sharp_6$ denotes the sum of the
number of 6-secants to $S$ and the number of $(-1)$-lines on $S$, then the
formulas
of Le Barz yield the following numbers for a surface of degree
10.\par\vfil\eject
$\pi =9$:

$$\matrix {&\quad \chi =1&\chi=2&\chi=3\cr
\sharp_5&6&12&18&\cr \sharp_6&7&3&3\cr}$$

$\pi =10$:$$\matrix {&\quad \chi =3&\chi=4\cr
\sharp_5&2&6\cr \sharp_6&2&1\cr}$$
  \bigskip

{\bf  Linkage  0.9.} ([PS])  Two surfaces $S$ and $S^{\prime}$ are said to
be linked $(m,n)$ if there exist hypersurfaces $V$ and $V^{\prime}$ of
degree $n$ and $m$ respectively such that $V\cap V^{\prime}=S\cup S^{\prime}$.
There are the standard sequences of linkage, namely $$\es SK{S\cup
S^{\prime}}{m+n-5}{S^{\prime}}{m+n-5}$$ $$\es SKS{m+n-5}{S\cap
S^{\prime}}{m+n-5}$$ The first sequence yields the relation between the
Euler-Poincar\'e characteristics $$\chi (S^{\prime})=\chi (V\cap
V^{\prime})-\chi (\kn S{m+n-5}). \leqno {(0.10)}$$The corresponding
sequence for linkage of curves in $\Pn 3$ yields the following relation between
the
sectional genera. $$\pi (S)-\pi (S^{\prime})={1\over
2}(m+n-4)(d(S)-d(S^{\prime})).\leqno { (0.11)}$$ To determine the
surfaces to which our surfaces are linked with, we will use:

\proclaim  Proposition 0.12.  If $S$ and $T$ are linked, then $S$ is
locally Cohen-Macaulay if and only if $T$ is locally Cohen-Macaulay.\par
$Proof.$ See [PS, Proposition 1.3].\par\smallskip

 For a proof of existence via linkage, the following
proposition will be used.

\proclaim Proposition 0.13.  If $T$ is a local complete intersection
surface in $\Pn 4$, which scheme\-theoretically is cut out by
hypersurfaces of degree $d$, then $T$ is linked to a smooth surface
$S$ in the complete intersection of two hypersurfaces of degree
$d$.\par
$Proof.$ See  [PS, Proposition 4.1].\par\bigskip

{\bf Remark} (Peskine, private communication).  A slight modification of
the conditions of this proposition is allowable, without changing the
conclusion. Namely, at a finite set of points $T$ need not be a local complete
intersection.  It suffices that it is locally Cohen-Macaulay, and that the
tangent cone at that point is linked to a plane in a complete intersection.
\smallskip
 The proof is an application of the proof of (0.13) to the strict transform of
$T$ in the blow-up of $\Pn 4$ in the points where $T$ is not a local complete
intersection.\medskip
Some useful lemmas:

\proclaim Lemma 0.14. Let  $E$ be a group of $t\leq 12$  points
 in $\Pn 2$, some possibly infinitely close, and assume that $\cohi 1E4\geq 1$,
where ${\cal I}_E$ is the sheaf of ideals defining the scheme $E$. Then either
$E$ is
a complete intersection of a cubic and a quartic curve, or there exists a
subgroup $E^{\prime}\subset E$ consisting of 6 points on a line, or of 10 or 11
points on a
conic.\par

$Proof.$ (cf. [EP, Cor 2 and Remark]).$\slutt$ \medskip
To check whether a linear system is very ample we'll use the following lemma,
which was communicated to us by J. Alexander.

\proclaim Lemma 0.15.  If $H$ has a decomposition
$$H\equiv C+D,$$ where $C$ and $D$ are curves on $S$, such that dim$|C|\ge 1$,
and if the restriction maps $\Coh 0SH\to\Coh 0DH$ and  $\Coh 0SH\to\Coh 0CH$
are surjective, and $|H|$ restricts to very ample linear systems on $D$ and on
every $C$ in $|C|$, then $|H|$ is very ample on $S$.\par
$Proof.$  We use the decomposition $H\equiv C+D$ to show that $|H|$ separates
points
and tangent directions on $S$.  Let $p$ and $q$ be two, possibly infinitely
close,
points on $S$.  By the assumptions of the lemma we may assume that $p+q$ is not
contained in $D$ or any $C$.  In particular we may assume that $p+q$ does not
meet
the base locus of $|C|$.  If $D$ contains $p$, then we can find a  curve $C$
which does not meet $p+q$ such that $C+D$ separates $p$ and $q$.  If $D$ does
not
meet $p+q$, then we can find a curve $C$ which contains one of the points $p$
or
$q$, such that  $C+D$ separates $p$ and $q$. $\slutt$\medskip
We deal with surfaces of degree 10 based on the following

\proclaim Proposition 0.16.  If $S$ is a smooth surface of degree 10 in $\Pn 4$
and
$\pi$ denotes the genus of a general hyperplane section,
then\smallskip\noindent
$\pi =6$ and $S$ is abelian or bielliptic, or\smallskip\noindent
$\pi =8$ and $S$ is an Enriques surface with four $(-1)$-lines, or a rational
surface, or\smallskip\noindent
$\pi =9$ and $S$ is a rational surface, or a blown-up $K3$
surface, or an honestly elliptic surface with $p_g=1$, $q=0$
and with three $(-1)$-lines, or a minimal surface of general type with $p_g=2$,
$q=0$, $K^2=3$ and one $(-2)$-curve, or\smallskip\noindent
$\pi=10$ and $S$ is a proper elliptic surface with $p_g=2$, $q=0$ and two
$(-1)$-lines,
or a minimal surface of  general type with $p_g=3$, $q=0$, $K^2=4$ and three
$(-2)$-curves,
or \smallskip\noindent
$\pi =11$ and $S$ is linked to an elliptic quintic scroll
($S$ lies on a cubic hypersurface), or $S$ is linked to a Bordiga surface
($S$ does not lie on a cubic hypersurface), or\smallskip\noindent
$\pi =12$ and $S$ is linked to a degenerate quadric surface,
or\smallskip\noindent
$\pi =16$ and $S$ is a complete intersection of a quadric and a quintic
hypersurface.\par
$Proof.$ (cf. [Ra]).$\slutt$\bigskip\eject

{\bf 1\quad 6-secants and plane curves}\bigskip

\proclaim  Lemma 1.1.  If $S$ contains a plane curve of degree $d_p$ and
$p_g=0$, then $\coh 1SH\geq {1\over 2}(d_p-2)(d_p-3)$, while if $p_g\geq 1$,
then $\coh 1SH\geq {1\over 2}(d_p-2)(d_p-3)+1-p_g$  .\par
$Proof.$  Let $C$ be a
plane curve on $S$, and let $D\in |H-C|$, and consider the cohomology of
the exact sequence $$\es SDSHCH.$$ If $p_g=0$, then $\coh 2SD=0$, so $\coh
1SH\leq \coh 1CH$. \par
If $p_g\geq 1$, then $\coh 2SD\leq p_g-1$, so $\coh 1SH\geq\coh
1CH-p_g+1$.$\slutt$

\proclaim  Corollary 1.2.  If $S$ has degree 10 and
sectional genus 9 or 10, then any plane curve on $S$ has degree at
most 4.\par
$Proof.$ Immediate, since these surfaces have $\coh 1SH \leq
2$.$\slutt$\par

\proclaim Lemma 1.3. If $S$ has a plane quartic curve $C$ with
$C^2=1$, then $\pi = 9$ and $S$ has three 6-secants in the plane of
$C$.\par
$Proof.$ The pencil of curves $|D|=|H-C|$ has $D^2=3$ base points in the
plane of $C$. The general member is smooth of genus 4 if $\pi =9$, and of genus
5 if
$\pi=10$.  The latter is impossible by the genus bound.  In the former case the
lemma follows, unless the base points of $|D|$ are collinear. Now the general
curve $D$ is a complete intersection of a cubic and a quadric. Therefore the
base points of $|D|$ are collinear if and only if the points $D\cap C$ are also
collinear.  But this means that $\coh 0C{H-D}>0$, while $\coh 0SC=1$ and $S$ is
regular by proposition 0.16, so this is impossible.$\slutt$

\proclaim Lemma 1.4. $\cohi 2Sn =0$ and $\cohi 2Hn=0$ for any hyperplane
section $H$ of $S$, when $n\geq 2$.\par

$Proof.$  If $\pi =9$ (resp. $\pi=10$), then  $1\leq\chi\leq 3$
(resp. $3\leq \chi \leq 4$), so $\cohi 3S1=0$ and
by Severi's theorem, $0\leq\cohi 2S1\leq 2$ (resp. $0\leq\cohi 2S1\leq 1$).
 Thus if $ \cohi 2S2 >0$ and $\chi\geq 2$, then $\cohi 2H2 >0$ for at least a
web
of hyperplane sections $H$.  If $ \cohi 2S2 >0$ and $\chi =1$ and $\pi =9$,
then $\cohi 2H2 >0$ for at least a net of hyperplane sections $H$.    In the
first case the general hyperplane section $H$ in the web is smooth, but $2\pi
-2<20$ so $\kn H{2H}$ is non special, i.e., $\coh 1H2=\cohi 2H2 =0$, a
contradiction.
In case $\chi =1$, $\pi =9$ and the general $H$ in the above net is not
smooth, i.e., the net has a line as a fixed curve, $H$ decomposes as
$H=L+C$, where $L$ is a line and $C$ is smooth. Consider the cohomology of the
exact
sequence
$$\es C{2H-L}H{2H}L{2H}.$$
$\coh 1H{2H}=\cohi 2H2 >0$, while $\coh 1L{2H}=0$, so
$\coh 1C{2H-L}>0$.  Therefore deg$\kn C{2H-L}=18-C\cdot L=17+L^2\leq 2p-2$,
where $p$ is the arithmetic genus of $C$. But $2p-2=2\pi -2 +2-2C\cdot
L=16+2L^2$, so  $17+L^2\leq 16+2L^2$, i.e., $L^2\geq 1$, which is impossible.
Again since $\cohi 3S1=0$ the second part of the lemma follows for $n=2$.
The argument applies inductively to the cases $n\geq 3$.$\slutt $\medskip

\proclaim Lemma 1.5. If $\pi =9$ or $\pi =10$, then $\cohi 0H3\leq 1$ for
any hyperplane section $H$ of $S$. Furthermore  $2\leq\cohi 1H3\leq 3$ when
$\pi=9$, while $1\leq\cohi 1H3\leq 2$ when $\pi=10$. \par
$Proof.$  If $H$ is contained in several
cubics, then the cubics are reducible, thus either $H$ decomposes into a line
$L$
and a curve $C$ of degree 9 which lies on a quadric with $L\cdot C\leq 2$, or
$H$
has a plane curve component of degree at least 6.  In the first case
$L^2=(H-C)\cdot L=1-L\cdot C\geq -1$ so $L^2=-1$ and $C$ has arithmetic
genus 8 (resp.9), which is impossible on a quadric, while the second case
contradicts lemma 1.1. The second part of the lemma now follows from
Riemann-Roch and lemma 1.4. $\slutt$\medskip

 Let  $\deh n$ be the locus in $\Pnd 4$ where
$(\cohi 0Hn)\cdot (\cohi 1Hn)\not= 0$.  This is clearly a determinantal
variety. To study these degeneracy loci we use the following observations:

\proclaim Lemma 1.6.  If $S$ has a $6$-secant, then $\deh 4$ contains a
plane.\par
$Proof.$ For the general plane $\Pi $
through the 6-secant $\cohi 1{\Pi\cap S}4>0$, while $\cohi 2H3=0$
for any hyperplane section containing $\Pi\cap S$ by lemma 1.4, whence
$\cohi 1H4>0$ for the general hyperplane $H$ through the 6-secant,
and the lemma follows. $\slutt$
\proclaim Lemma 1.7.  If $H\in \deh 4$, then $H$
has a proper $6$-secant or a plane curve component. \par
$Proof.$ If $H$ does not have a plane curve component, then the sequence
$$\esi H3H4{\Pi\cap H}4$$ is exact for any plane $\Pi$ in this hyperplane.
Since $\cohi
1H3\leq 3$, there exists at least one plane section for which $\cohi 1{\Pi\cap
S}4>0$.  This plane section is a scheme of length 10. If it is contained in
a conic, then $\cohi 1H1>0$, but $S$ is regular so this means that $\cohi
1S1>0$, which contradicts  Severi's theorem.  Therefore it follows from lemma
0.14
that $\Pi\cap S$ contains a subscheme of length 6 which is contained in a line,
and the lemma follows. $\slutt$
\proclaim Lemma 1.8.   $S$ has only finitely many plane
curves.\par
$Proof$. Since $S$ is regular, any one-dimensional family of plane
curves on $S$ is linear. Thus the curve
residual to the general plane curve in a hyperplane is again a plane curve,
so the hyperplane section is contained in a quadric, impossible by Severi.
$\slutt$
\proclaim Lemma 1.9. Every component of $\deh 4$ of dimension at least $1$
contains a line through each of its points.  Every component of $\deh 4$ of
dimension at least $2$ contains a plane through each of its points. A
plane in $\deh 4$ correspond to a proper 6-secant or to a line contained in
$S$.\par
$Proof.$  By lemma 1.6, lemma 1.7 and lemma 1.8 every point in a
component of dimension at least 1 lies in a plane or on a line whose
corresponding hyperplane sections have a common plane component. Any
component of dimension at least 2 without planes through each point
gives rise to a family of plane curves, thus contradicting lemma 1.8. If the
hyperplane sections corresponding to a plane have no line in common,
then the general one is smooth and has a
proper 6-secant by lemma 1.7. $\slutt$
\proclaim Lemma 1.10. $\deh 4$ cannot contain a $\Pn 3$ or  a determinantal
quadric or cubic hypersurface.\par
$Proof$. In either case $\deh 4$ contains at least a pencil of pairs of
planes which meet along a line. Since $S$ is not a scroll, it must contain
infinitely many 6-secant lines by lemma 1.9. Furthermore, a pair of 6-secants
in the pencil
meet, since the corresponding planes in $\Pnd 4$ meet in a line. By Bezout this
means that $S$ meets the plane spanned by a general pair of 6-secants in a
curve.
 Thus there is a pencil of plane curves on $S$, impossible by lemma 1.8.$\slutt
$

\proclaim Proposition 1.11.  Let $S$ be a smooth surface of degree $10$ and
sectional genus $9$.  If $\chi =1$, then the number $N_6$ of proper $6$-secants
is at most $1$, and  $\cohi 1S4$ is also $0$ or $1$. If $\chi =2$, then
$N_6=0$ or $N_6=1$ and $\cohi 1S4=1$, or $N_6=3$ and  $\cohi 1S4=2$.
If $\chi =3$, then $N_6=3$ and $\cohi 1S4=2$.  When the number of 6-secants
is $3$, then $S$ cuts some plane along a quartic curve with self-intersection 1
and three
points.\par
\proclaim Remark 1.12.  The formula of Le Barz gives the sum of the number
of $6$-secants and the number of $(-1)$-lines on $S$.  When $\chi =1$ the sum
is
$7$, and when $\chi =2$ and $\chi=3$ the sum is $3$.\par

  First note that $\deh 3\subset \deh 4$, since if $H$ lies
on a cubic, then it lies on at least 4 quartics, which means that $\cohi
1H4>0$. Let $\deh {3,4}$
denote the loci of hyperplanes for which $\cohi 1H3 =3$ and $\cohi 1H4=2$,
or equivalently $\cohi 0H3=1$ and $\cohi 0H4=5$.
Now $\cohi 1H3 =2$ and $\cohi 1H4\geq 2$, or $\cohi 1H3 =3$ and $\cohi 1H4\geq
3$
are impossible like in lemma 1.10.  Thus
\proclaim Remark 1.13. $H\in \deh
{3,4}$ if and only if $\cohi 1H4=2$.\par
\proclaim Lemma 1.14. If $\cohi 1H4=2$,
then $H$ is reducible; it has a plane quartic component $C$ and a connected
residual component $D$, which is a complete intersection (2,3) or decomposes
into a
plane quartic component (which may coincide with $C$) and a conic.
Furthermore,
$S$ has three 6-secants in a plane.\par
$Proof.$  First note that $H$ is contained
in a cubic and an independent quartic hypersurface.
 If $H$ is linked $(3,4)$ to a
curve $C$ of degree 2 and arithmetic genus -3, then the curve $C$ must be a
double structure on a line.  But any such structure on a cubic has
arithmetic genus at least -2, a contradiction. So the cubic and the
quartic containing $H$ must have a common component. So $H$ has a plane
component, of degree at least 2 and at most 4, and possibly a component on
an irreducible quadric.  Any such component $C$ on a quadric is of type (a,b)
with $a\leq b\leq 5$, $a+b\leq 8$; for if say $b=6$ and $a=2$, then $C$ has
arithmetic genus $p_a(C)=5$ and since the sectional genus is 9, $C$ must meet
the residual curve $H-C$ of degree 2 in $C\cdot H-C\geq 5$ points, which means
that $H-C$ has a component on the quadric, the other values of $a$ and $b$
with $b\geq 6$ are impossible by similar reasoning.\par
 If the cubic and
the quartic have a common irreducible quadric, then the plane curve is a
conic A, and the curve $B$ on the quadric is of type (4,4) or (3,5), but
$A^2\leq -1$ so $A\cdot B\geq 3$, which means that $p_a(A+B)\geq 10$,
impossible.
Therefore the quartic and the cubic have a plane in common. The curve $C$
of this plane is residual to a curve $D$ in $H$ which is on a quadric and an
independent cubic. The quadric and the cubic have at most a plane in common,
and any
plane curve has degree at most 4 by corollary 1.2, so the curve $D$ has
degree at most 6, and thus $C$ is a plane quartic curve.   The curve $D$ moves
in a pencil on $S$, whose general element is irreducible, so $D$ is connected.
If
the quadric and the cubic have no common component, then $D$ is a
complete intersection (2,3). If the quadric and the cubic has a common
plane, then the curve $A$  of this plane must be a quartic, since the
residual curve $B$ is a conic lying in a plane and an independent quadric.  Now
$1\leq B\cdot A\leq 2$ and $C\cdot A\leq 2$, while  $B^2\leq -1$
implies that $B\cdot (A+C)\geq 3$.
  If $B\cdot A= 2$ or $D$ is a complete intersection, then $p_a(D)=4$ and
$C\cdot D=3$ and  $C^2=H\cdot C-C\cdot
D=1$ so, by  lemma 1.3,  $S$ has three 6-secants in the plane of
$C$.  If $B\cdot A= 1$, then $C\cdot A=2$ and $B+C$ has genus 4.
As above $A^2=1$ and there are three 6-secants in the plane of the quartic $A$.
$\slutt$\medskip

\proclaim Lemma 1.15.  $\cohi 0S3=0$ and $\cohi
1S3=\chi+1$.\par
$Proof.$ Now $\cohi 0S3\leq 1$ by lemma 1.5.  If $\cohi 0S3=1$, then $\cohi
1S3=\chi+2$, while  $\cohi 1S2=\chi -1$.  So if $\chi =2$ or $\chi =3$,
then $\cohi 0H3=2$ for some $H$, contradicting lemma 1.5.  If
$\chi =1$, then $\cohi 1S3= 3$ and $\cohi 1S4\geq 4$, so $\cohi 1H4=2$ for
at least a 2-dimensional
 family of hyperplane sections. By lemma 1.14, $S$ has at least a
1-dimensional family of plane quartic curves, which is impossible by lemma
1.8. $\slutt$ \medskip

$Proof$ $of$ $proposition$ $1.11.$   When $\chi =1$, then $\cohi 1S3=2$.
So if $\cohi 1S4>1$, then $\deh 4 $ contains a $\Pn 3$ or a rank 4 quadric,
thus contradicting lemma 1.10.
When $\chi =2$,  then $\cohi 1S2= 1$ and $\cohi 1S3= 3$ by lemma 1.15,
so  $\cohi 1H3= 3$ and $\cohi 1H4\geq 1$ for some $H$.
Therefore $\cohi 1S4\geq 1$.   If $\cohi 1S4=1$, then $\cohi 1H4\leq 1$,
and $\deh 4$ is defined by three linear forms. So $\deh 4$ is a
line
or a plane.  If it is a line, $S$ has no 6-secants. If $\deh 4$ is a plane,
then $S$ has at most one 6-secant.
 If $\cohi 1S4=2$, then $\deh 4$ is defined by the minors of a $2\times 3$
matrix with linear entries.  Thus $\deh 4$ is the union of three planes by
lemma 1.9 and lemma 1.10,  and  $\cohi 1H4=2$ for
some hyperplane section $H$.  By lemma 1.14, $S$ has three 6-secants in the
plane
of a plane quartic and a residual pencil of curves of degree 6 and genus 4.
If $\cohi 1S4\geq 3$, then $S$ has too many plane curves. Now if $\chi
=3$, then $\cohi 1S2= 2$ and $\cohi 1S3= 4$ by lemma 1.15.  Recall that $\deh
3\subset
\deh 4$, hence $\cohi 1S4\geq 2$. Since $S$ is minimal, it follows from
remark 1.12. that $S$ has 3 or infinitely many 6-secants.  Therefore $\deh 4$
contains
at least three planes, so if $\cohi 1S4=2$, then $\deh {3,4}$ is not empty,
and lemma 1.14 applies, to show that $S$ has three 6-secants in a plane.   If
$\cohi 1S4\geq 3$, then $\deh 4$ contains a determinantal cubic, thus
contradicting lemma
1.10.$\slutt $

\bigskip
\proclaim Proposition 1.16.  Let $S$ be a smooth surface of degree 10 and
sectional genus 10. Then
$$(\cohi 0Sn)\cdot (\cohi 1Sn)=0$$ when $0\leq n\leq 4$ and $\chi =3$,
while
$$(\cohi 0Sn)\cdot (\cohi 1Sn)=0$$ when $0\leq n\leq 3$, and $\cohi 0S4 =4$,
$\cohi 1S4=1$ when $\chi =4$.  In particular, $S$ has a 6-secant when $\chi
=4$, and it has no 6-secant when $\chi =3$.\par \proclaim Remark 1.17.
The formula of Le
Barz gives the sum of the number of 6-secants and the number of
$(-1)$-lines on $S$.  When $\chi =3$ the sum is 2, and when $\chi =4$
the sum is 1.\par
\proclaim Lemma 1.18.  $\cohi 0S3=0$ and $\cohi
1S3=\chi -2$\par
$Proof.$  $\cohi 0S3\leq 1$ by lemma 1.5.
  If $\cohi 0S3=1$, then, by lemma 1.4, $\cohi 1S3= \chi
-1$, while $\cohi 1S2=\chi -3$.  So if $\chi =4$, then $\cohi
1H3=3$ for some $H$, contradicting lemma 1.5.  Furthermore, $\cohi 0S3=1$
implies
that $\cohi 0S4\geq 5$, so if $\chi = 3$, then $\cohi 1S4 \geq 2$.
Therefore $\deh 4$ contains a determinantal quadric, thus contradicting lemma
1.10. $\slutt$\medskip
$Proof$ $of$ $proposition$ 1.16. Now $\cohi
1S3=\chi -2$, by lemma 1.4 and lemma 1.18.  If $\chi =3	$, then $\cohi
1S4=0$, by lemma 1.10. If $\chi =4$, then $\cohi 1S2=1$ and $\cohi 1S3=2$ so
$\deh
3$ is at least a plane.  Recall that $\deh 3\subset \deh 4$, so $\cohi 1S4\geq
1$.  But $\cohi 1S4\leq 1$, by lemma 1.10.  Therefore $\cohi 1S4=1$ and $\deh
4$
is a plane.  By remark 1.17, $S$ has one 6-secant since $S$ is minimal when
$\chi = 4$. When $\chi =3$, $\deh 4$ is empty, so
$S$ has no 6-secants in this case.$\slutt$\bigskip\bigskip

{\bf 2\quad Linear systems}\bigskip

\proclaim Proposition 2.1. If $S$ is a smooth
rational surface of degree 10 in $\Pn 4$ with $\pi =9$, then
$$H\equiv8\pi^{\ast}\li -\sum_{i=1}^{12}2E_i-\sum_{j=13}^{18}E_j,$$
 or  $$H\equiv 9\pi^{\ast}\li
-\sum_{i=1}^43E_i-\sum_{j=5}^{11}2E_j-\sum_{k=12}^{18}E_k,$$
where $\pi:S\to \Pn 2$ is the blow-up map with exceptional
curves $E_i$, $i=\overline {1,18}$, and $\li $ is a line in $\Pn
2$.\par\bigskip

$Proof.$  From proposition 1.11 and the 6-secant formula
(0.8) it follows that the number of $(-1)$-lines on $S$ is 6 or 7.   The
procedure is first to produce a number of possible candidates using
adjunction,
and secondly to show that only those described in the proposition are
possible.\medskip  For $\pi =9$ we obtain the following list of invariants for
$S$, $S_1$ and $\Sigma$, where $S_1$ is the image of $S$ under the
adjunction map and $\Sigma$
is the image of $S_1$ under the adjunction map defined by
$|H_1+K_1|$:\medskip
\settabs \+$S\subset \Pn 4$\qquad&$  H^2=10$\qquad\qquad&$  H\cdot
K=6$\qquad\qquad\qquad&$ K^2=-9$\qquad\qquad\qquad&$\pi =9$\cr\+$S\subset
\Pn 4$&$  H^2=10$&$  H\cdot K=6$&$ K^2=-9$&$\pi =9$\cr\+$S_1\subset\Pn 8$&$
H_1^2=13$&$  H_1\cdot K_1=-3$&$ K_1^2=-9+a$&$\pi_1=6$\cr\+ $\Sigma\subset
\Pn 5$&$  H_{\Sigma}^2=a-2$&$  H_{\Sigma}\cdot K_{\Sigma}=a-12$&$
K_{\Sigma}^2=-9+a+b$&$\pi_{\Sigma}=a-6,$\cr  \medskip\noindent where $a$ is
the number of $(-1)$-lines on $S$ and $b$ is the number of $(-1)$-lines on
$S_1$.
$\Sigma $ is a surface since $4\leq H_{\Sigma}^2=a-2\leq 5$ when $6\leq a\leq
7$.
 So the invariants for $\Sigma$ make sense.  When $H_{\Sigma}^2=4$, the
surface $\Sigma$ is a Veronese surface or a rational normal scroll, and when
$H_{\Sigma}^2=5$ the surface $\Sigma$ is a Del Pezzo surface. Thus $H$ can
be reconstructed via the adjunction process to get the following list of
candidates:\par\bigskip
\hskip 2.3cm$\eqalign{ 1)\quad\quad\quad H&\equiv
5B+(6-{5\over2}e)F-\sum_{i=1}^{11}2E_i-\sum_{j=12}^{17}E_j,\qquad e=0\;{\rm
or}\;2,\cr}$\par\bigskip
where $B$ is a section with self-intersection $B^2=e$ and $F$ is a member of
the ruling,\par\bigskip
\hskip 2.3cm $\eqalign{2)\quad\quad\quad H&\equiv8\pi^{\ast}\li
-\sum_{i=1}^{12}2E_i-\sum_{j=13}^{18}E_j,\cr 3)\quad\quad\quad H&
\equiv 9\pi^{\ast}\li
-\sum_{i=1}^43E_i-\sum_{j=5}^{11}2E_j-\sum_{k=12}^{18}E_k.\cr }$\par\bigskip
In case 1) we study curves in the linear system
$$|C|=|2B+(3-e)F-\sum_{i=1}^{11}E_i|.$$ Since dim$|2B+(3-e)F|$=11, there is
a curve $C$ in $|C|$. It has degree 5 and arithmetic genus 2.  If $C$ is not
contained in a hyperplane, then it must be the union of a plane quartic $A$
and a line $L$ not meeting the plane of $A$. By the index theorem $A^2\le 1$,
so $L^2=1-A^2\ge 0$, which is impossible by Riemann-Roch.  Therefore $C$ is
contained in a hyperplane section $H$ with a residual curve $C_1=H-C$ of
degree 5 and arithmetic genus 4, which is
impossible by lemma 0.3.$\slutt$\bigskip\noindent

 \proclaim Proposition 2.2.  If $$|H|=|8\pi^{\ast}\li
-\sum_{i=1}^{12}2E_i-\sum_{j=1}^{6}F_j|$$ is the linear system of
hyperplane sections of a smooth surface $S$ in $\Pn 4$, then there are three
plane
quartic curves on S, whose respective planes all contain the same line, which
in
turn is a 6-secant line for the surface.  The map $\phi
C$ defined by $$|C|=|4\pi^{\ast}\li -\sum_{i=1}^{12}E_i|,$$
is of degree 4 onto $\Pn 2$, and maps a member of the linear system
$|\pi^{\ast}\li|$ to a plane quartic with three nodes such that each node is
the image by $\phi C$ of two of the exceptional curves $F_j$.\par

\medskip\noindent
$Proof.$  Let $$C\equiv 4\pi^{\ast}\li -\sum_{i=1}^{12}E_i,$$ $$
C_{ij}\equiv C-\sum_{k=1}^{6}F_k+F_i+F_j$$
  and
$$C^{ij}\equiv H-C_{ij},\qquad {\rm for}\quad 1\le i<j\le 6.$$  Now $\coh
0S{C^{ij}}>0$, so there is a
curve $C^{ij}$ in $|C^{ij}|$.  It has degree 6 and arithmetic genus 3, so
by lemma $0.3$, it either spans a $\Pn 3$, in which case there is a residual
plane
quartic curve $C_{ij}\equiv H-C^{ij}$ and $\coh 0S{C^{ij}}=2$, or $C^{ij}$ is
reducible, that is $C^{ij}$ is the union of a plane quartic $A$ and two skew
lines
$L_1$ and $L_2$, with $A\cdot L_1=A\cdot L_2=1$, or a conic $Q$, with $A\cdot
Q=1$.
One may now check that in the latter case $A\equiv C_{st}$ for some $s,t$ and
$L_1$ and $L_2$ are lines $F_k$ and $F_l$.\medskip\noindent  These two
possibilities for each $C^{ij}$ fit together only if say  $$C_{12},\quad C_{34}
\quad{\rm and}\quad C_{56}$$ are plane quartics and  $$|C^{12}|,\quad
|C^{23}|\quad{\rm and}\quad |C^{56}|$$ are their respective residual
pencils.\medskip\noindent
Now $C_{12}\cdot C_{34}=C_{12}\cdot C_{56}=C_{34}\cdot C_{56}=2$, so the
planes of
$C_{12},\;\;C_{34}$ and $C_{56}$ meet pairwise in lines.  Since the three
planes span all of $\Pn
4$, they must intersect in a common line $L$, which is now a 6-secant for
the surface $S$ unless $L$ lies on $S$.
\medskip\noindent To see that $L$ cannot lie on $S$, we first note that $S$
has no plane quintic curve by $(1.2)$.  So if $L$ lies on $S$, then $L$ is a
component of the curves $C_{ij}$ and $L\cdot (C_{ij}-L)=3$.  In this case
we get
$C_{12}\cdot C_{34}=(C_{12}-L)\cdot (C_{34}-L)+6+L^2=2$, but
$(C_{12}-L)\cdot
(C_{34}-L)\geq 0$ so $L^2\leq -4$.  Thus $C_{12}\cdot L=C_{34}\cdot L\leq
-1$.    Since $H\cdot L=(C_{12}+C_{34}+F_5+F_6)\cdot L=1$, this means that
$L\cdot F_5>1$ or $L\cdot F_6>1$, which is absurd. So $L$ cannot be contained
in
$S$.

 \medskip\noindent
Let $S\cap L= q_1+\dots+q_6$ be six, some possibly infinitely
close, points such that $C_{12}\cap L=q_3+q_4+q_5+q_6$, $C_{34}\cap
L=q_1+q_2+q_5+q_6$ and $C_{56}\cap
L=q_1+q_2+q_3+q_4$. Since $C_{12}$ is a plane quartic curve on $S$, that is
$$\kn {C_{12}}H\cong \omega_{C_{12}}\cong \kn {C_{12}}{\pi^{\ast}\li},$$ we see
that the collinear points $q_3,\dots,q_6$ all lie on a curve
$L_0\equiv \pi^{\ast}\li$.
But we get the same one for $C_{34}$ and $C_{56}$, so this means that all the
points $q_i$ lie on $L_0$.  The map $\phi C$ maps the base points of a
pencil in $|C|$ to a point, so $q_1,q_2$ and $F_1,F_2$ are mapped to the
same point which is a node of $\phi C(L_0)$, etc.   $\slutt$\medskip
The quickest construction of this surface is by linkage, cf. (4.3).

\proclaim Proposition 2.3.  Assume that $$|H|=|9\pi^{\ast}\li
-\sum_{i=1}^{4}3E_i-\sum_{j=5}^{11}2E_j-\sum_{k=12}^{18}E_k|$$ is the
linear system of hyperplane sections of a smooth surface $S$ in $\Pn 4$.
Let $\Sigma$ be the image of $S\subset \Pn 3$ under the map $\phi C$
defined by $$|C|=|4\pi^{\ast}\li -\sum_{i=1}^{11}E_i|.$$  There is
a twisted cubic curve in this $\Pn 3$ for which the lines $\phi C(E_i)$,
$i=\overline{1,4}$, are secants and which meets $\Sigma $ in the
points $\phi C (E_k)$, $k=\overline{12,18}$, outside these
lines. Then $|H|$ is given by the linear system of quartic surfaces through
the double curve of $\Sigma$, the four lines $\phi C(E_i)$,
$i=\overline{1,4}$, and the seven points $\phi C(E_k)$,
$k=\overline{12,18}$. Moreover, $S$ lies on a determinantal
quartic with 36 nodes, all lying on $S$.\par
$Proof.$ First a few lemmas.
\proclaim Lemma 2.4.  The linear
system $|C|=|4\pi^{\ast}\li -\sum_{i=1}^{11}E_i|$ has dimension 3 and has
no base points, the double curve of its image is a, possibly reducible, twisted
cubic curve and it has no triple points.\par
$Proof$.  If dim$|C|>3$,
then, by lemma $0.14$, the linear system has a fixed curve which will have
negative degree on $S$, absurd, so dim$|C|=3$.  Similarly, any fixed curve
in the system will again have negative degree. If it has a basepoint, then
again
lemma $0.14$ applies to show that there is an elliptic curve $3\pi^{\ast}\li
-\sum_{i=1}^{11}E_i$ of degree one on $S$, which is also absurd.  So the
first part of the lemma follows.  Any smooth member of the system is a
non-hyperelliptic curve of genus 3 which is mapped to a plane quintic curve
by $\phi C$.  Since the curve is not hyperelliptic, the image has only double
points as singularities.  Thus the double point class $7\pi^{\ast}\li
-\sum_{i=1}^{11}2E_i$ contains a curve $D_C$ on $S$ which is mapped
two to one onto its image. And the image is a curve of degree 3 which spans
$\Pn 3$.  It is residual to a rational quartic curve in the intersection of
$\phi C(S)$ and a quadric surface in $\Pn 3$. Thus it is a twisted cubic
curve.$\slutt$

\proclaim Lemma 2.5.  The linear
system $|H|$ is cut out on $\Sigma$ by quartic surfaces containing
the twisted cubic curve  $\phi C(D_C)$, the four secants $\phi C(E_i)$,
$i=\overline{1,4}$, and the points $\phi C(E_{k})$, $k=\overline{12,18}$,
i.e., $|H|$ extends to a linear system on a blowup of $\Pn 3$.\par

$Proof$. Consider the linear system of quartic
surfaces in $\Pn 3$ which contains the twisted cubic curve  $C_3=\phi C (D_C)$
and the four secants $\phi C (E_i)$ $i=\overline {1,4}$.  Denote by ${\cal M}$
the
moving part of the pullback of this system to $S$.  Clearly $${\cal M}\subset
|H_0|=|9\pi^{\ast}\li -\sum_{i=1}^{4}3E_i-\sum_{j=5}^{11}2E_j|.$$
Clearly the lemma follows if this is an equality. But the union of a
twisted cubic curve and four of its secants is a
curve of arithmetic genus at least 4, so it is contained in at least 10
quartic surfaces.  Therefore dim${\cal M}\geq 9$.  On the other hand, $S$ is
regular and $\coh 0{H_0}{H_0}=9$ for any smooth curve $H_0$,
so dim$|H_0|=9$.$\slutt$\medskip

We proceed to show that the base locus of this linear system of quartic
surfaces is an arithmetically Cohen Macaulay curve of genus 11 and degree 10.
\proclaim Lemma 2.6.  The general member $D$ of the linear system
$|12\pi^{\ast}\li-\sum_{i=1}^{4}4E_i-\sum_{j=5}^{11}3E_j-\sum_{k=12}^{18}E_k|$
is irreducible.\par

$Proof.$  Note that $D$ has degree 11 and arithmetic genus 10 and self
intersection $D^2=10$.  By Riemann-Roch $|D|$ is at least a pencil, so
the general $D$ is irreducible unless $|D|$ has a fixed curve. Let $F$
denote this fixed curve and let $M=D-F$ be the moving part of $D$.  By
Riemann-Roch $|2H-D|$ is also a pencil, so let us consider the space
$P$ of quadrics in $\Pn 4$ which contain $F$.  Since two members of
$D$ cannot be contained in the same quadric, the projective dimension
of $P$ is at least 3.  This reduces the number
of possibilities for $F$.\par
  First let us check whether $F$ can be contained in a hyperplane.
If so, we get the subcases a) $F$ is not on a quadric in
$\Pn 3$, b) $F$ is on a quadric, c) $F$ lies in a plane. In case a) the
quadrics containing $F$ decompose, so $M$ lies in a pencil of hyperplanes,
i.e., is a plane curve, but this means that $S$ has a pencil of plane
curves, contradicting lemma 1.8.   In case
b) there is at least a $\Pn 2\subset P$ of reducible quadrics, so every $M$
is also contained in a quadric in a $\Pn 3$, but since $M$ moves in a pencil,
this means that $A=H-M$ is a plane curve.
By lemma $1.2$ the degree of $A$ is at most
4. If it is 4, then $0\leq M^2\leq 1$ since $M$ moves and $S$ has no 6-secant.
But then, by (0.2) the curve $M$ has genus 1 or 2 and degree 6 on a quadric,
impossible.  If $A$ has degree 3, then $M$ must be of type (3,4) or (2,5)
on the quadric, while $F$ has degree 4 and genus at most 1.
 Now $M$ and $F$ meet in a
conic section, and by (0.2) we have $M\cdot F=11-p_a(M)-p_a(F)\geq 4$, with
equality if
$M$ is of type (3,4) and $F$ is elliptic.  But inequality means that $F$
has a component on the quadric of $M$, impossible.  With equality
$F^2=D^2-M^2-2M\cdot
F=-2$, since by (0.2) $M^2=4$. But then $H-F$ has degree 6 and genus 3, and
$(H-F)\cdot A=6$, impossible. If  $A$ is a conic section, then $M$ is of
type (3,5) or (4,4), which are both impossible by (0.2).  Similarly $A$ cannot
be a line. In case c) a similar argument to the one of case b) leads to a
contradiction.\par If $F$ is not contained in a hyperplane, then the
quadrics containing it cut out a curve, and $F$
has degree 6 and genus 2, or degree 5
and genus 0 or 1, or degree 4 and genus at most 0, or degree 3 and genus
at most $-1$.  In either case we can write $M=\alpha\pi^{\ast}\li
-\sum_{i=1}^{4}c_iE_i-\sum_{j=5}^{11}b_jE_j-\sum_{k=12}^{18}a_kE_k $.  Note
that $a_k$ is 0 or 1 for each i.  Furthermore, if $a_{12}+...+a_{18}<7$, then
one
member of
the pencil $|M|$ decomposes into a line and a curve $M_0$, where $M$ and
$M_0$ have the same arithmetic genus.    Thus $$M\cdot F=M_0\cdot
F+1\leq H\cdot F+1,$$ where $$M\cdot F=11-p_a(M)-p_a(F)$$ from (0.2).
Furthermore, if $M_2$ is the image of $M$ by the second adjunction, i.e.,
on the Del Pezzo, then
$$H_2\cdot M_2=-MK+\sum_{j=5}^{11}b_j+\sum_{k=12}^{18}a_k=H\cdot
M+4p_a(M)-4-2M^2-\sum_{k=12}^{18}a_k.$$\par
We go now case by case. If $p_a(F)=2$, then $M$ has degree 5 and genus at
most 2 since it is not a plane curve, so $H_2\cdot M_2\leq
5+4p_a(M)-4-\sum_{k=12}^{18}a_k\leq 5+4p_a(M)-4$ means that $M$ has genus 1
or 2.  If $p_a(M)=2$, then $M$ spans a hyperplane so $F\cdot M=7\leq F\cdot
H=6$, impossible. If $p_a(M)=1$, then $a_{12}+...+a_{18}<7$, so there is an
elliptic
curve $M_0$ of degree 4 which spans a hyperplane, and by the above $F\cdot
M_0=7\leq F\cdot H=6$, impossible. This argument works similarly when
$p_a(F)\leq
1$.$\slutt$

\proclaim Lemma 2.7. The
general member $D$ of the linear system $|12\pi^{\ast}\li
-\sum_{i=1}^{4}4E_i-\sum_{j=5}^{11}3E_j-\sum_{k=12}^{18}E_k|$ lies on a Del
Pezzo surface $S_D$ and is a member of the linear system $|7\pi^{\ast}\li
-\sum_{i=1}^{5}2E_i|$ on the Del Pezzo.\par $Proof.$ $D$ has degree 11 and
arithmetic genus 10.  By Riemann-Roch $D$ lies on at least two quadrics.
If it
lies on three quadrics, then it must lie on a cubic scroll. Since it lies
on an
independent quintic, it is linked to four skew lines on the scroll,  so it
is of
type $6\pi^{\ast}\li -e$ on the scroll.  The map to $\Pn 3$ is a $g^3_{11}$
which is not a projection of the map into $\Pn 4$.  So it is defined by a
system
$3\pi^{\ast}\li -\sum_{i=1}^7p_i$ on the cubic scroll, where the $p_i's $ are
assigned base points.  To have dimension 3 at least 5 of the $p_i's$ lie on a
line,
which means that the image in $\Pn 3$ lies on a quadric, absurd.  So $D$
lies on
two quadrics which defines a cone over an elliptic curve or a Del Pezzo.  In
either case it lies on an independent cubic, so it is linked to a line in
the
complete intersection of the quadrics and the cubic.  On a cone over an
elliptic
curve, the curve would have three branches through the vertex, imposing a
singularity at this point of the surface $S$.   So $D$ lies on a Del Pezzo,
and
is linked to a line in the complete intersection of the Del Pezzo surface
and a
cubic hypersurface.$\slutt$\par

\proclaim Lemma 2.8. The linear system $|C_{|D}|$ extends to a linear system
$|3\pi^{\ast}\li -\sum_{i=1}^{4}E_i|$ on
$S_D$ with exactly two base points on $D$.\par
$Proof.$  The adjoint linear
system to $D$ on $S_D$ is $|4\pi^{\ast}\li -\sum_{i=1}^{5}E_i|$. Let
$p_1+...+p_7$ be a general member of the pencil residual to $|C_{|D}|$ in
$|K_D|$, blow up $S_D$ in these points and let $F_i$, $i=\overline {1,7}$ be
the
exceptional curves.  Then $|4\pi^{\ast}\li
-\sum_{i=1}^{5}E_i-\sum_{i=1}^{7}F_i|$ has dimension 3.  By lemma 0.14 this
is
possible only if there is a curve $3\pi^{\ast}\li
-\sum_{i=1}^{5}E_i-\sum_{i=1}^{7}F_i$, a curve $2\pi^{\ast}\li
-\sum_{i=1}^{j}E_{k_i}-\sum_{i=1}^{10-j}F_{k_i}$ or a curve $\pi^{\ast}\li
-\sum_{i=1}^{j}E_{k_i}-\sum_{i=1}^{6-j}F_{k_i}$ on the blown up surface.
The
first case is impossible since the map ${\phi C}_{|D}$ would have a
quadruple
point, the second is impossible since then the image of $D$ would
lie on a quadric.  In the third case the fixed curve  must meet one of the
$E_i's$,
otherwise it would intersect $D$ and impose an unassigned basepoint,
impossible.  Thus the fixed curve  meets one of the $E_i's$, say $E_5$ and
the
lemma follows.$\slutt$\medskip

It follows immediately from this lemma that $\phi
C(D)$ lies on a cubic and has arithmetic genus 11.  The pencil of these
cubic surfaces, coming from the pencil $|D|$ has a base curve which consists of
the four lines $\phi C(E_i)$ $i=\overline {1,4}$, and their two transversals,
and a
twisted cubic curve for which the four lines are secants. Together with the
twisted cubic which is the double curve of $\Sigma =\phi C(S)$, the four lines
and this twisted cubic curve form a curve of degree 10 and genus 11, which
does not lie on a cubic (the cubic would belong to the above pencil).
Therefore it is defined by the maximal  minors of $4\times 5$-matrix
with linear entries. These minors restrict to $\phi C(S)$ to give
the embedding of $S$ in $\Pn 4$, i.e., the points $\phi C(F_i)$,
$i=\overline{1,7}$, are the intersection of $\phi C(S)$ with a
twisted cubic outside the lines  $\phi C(E_i)$, $i=\overline{1,4}$. The image
of $\Pn 3$ by this map is a determinantal quartic with
nodes corresponding to the 4-secants
and the singular points of the base curve of degree 10 and genus 11: Clearly
a 4-secant is contracted, it is a $(-2)$-curve on a general quartic in
the linear system containing it, so the image is a quadratic singularity. For
the singularities on the base curve we have

\proclaim Lemma 2.9. Let $C_1$ and $C_2$ be two curves intersecting
transversally in a point $p$ in $\Pn 3$.  Let $V\to \Pn 3$ be the blowup
along the ideal of $C_1\cup C_2$. Then $V$ has a quadratic singularity
at a point $q$ over $p$. The strict transform on $V$ of any surface
with multiplicity $a$ along $C_1$ and multiplicity $b$ along $C_2$
is locally Cohen-Macaulay at $q$ if and only if $|a-b|\leq 1$.
If the branches have normal crossings and are transverse to the plane
spanned by the two branches, then the strict transform is smooth at $q$ if
$|a-b|= 1$, and does not meet $q$ if $a=b$. \par

$Proof.$ This is a purely local calculation which can be worked out by blowing
up two concurrent lines in $\Pn 3$ and looking at the singular quadric in
$\Pn 4$ defined by the quadrics through the two lines. The details are left to
the reader.$\slutt$\par\medskip

 Now, the number of 4-secants is 20 and each of them intersect
$\Sigma$ once outside the base curve. The number of singular points is
16, and at each point the difference between the number of branches of
$\Sigma$ through the two branches of the base curve is 1, so the last
statement of proposition 2.3 follows.$\slutt$\bigskip

{\bf 2.10. Construction of B.}  To prove the existence it
is clear how to go the other way. Start with 11 general points in the
plane, or with a general quintic with a twisted cubic double curve. In both
cases the surface is smooth outside the twisted cubic curve and each point
on the twisted cubic is a double point of the surface.  Starting with the
plane the result follows from a study of the linear system $|4\pi^{\ast}\li
-\sum_{i=1}^{11}E_i|$.  Starting with four lines in $\Pn 3$ not all on a
quadric, pick two disjoint general pencils in the linear system of cubics
through these lines. Each pencil defines a twisted cubic curve
residual to the four lines and its two transversals. There is a pencil
of quartic surfaces double along the twisted cubic, and containing the
four lines, in fact the complete intersection is exactly the four
lines plus a quadruple structure on the twisted cubic curve. Thus it
is easy to see that the general quintic $\Sigma$, that we are looking
for, is smooth outside the twisted cubic curve.

Clearly the net of quadrics containing the twisted cubic defines a rational
map of degree one to $\Pn 2$, and the linear system on $\Pn 2$ can easily be
deduced from this map.   The union of the two twisted cubics and the four lines
form a curve of degree 10 and genus 11, which is defined by the maximal minors
of
$4\times 5$-matrix with linear entries.  The five maximal minors of this matrix
define a rational map of $\Pn 3$ into $\Pn 4$, the image is a determinantal
quartic.
 The following lemma gives a description of the singularities of this map.

  \proclaim Lemma 2.11.  Let $C\subset \Pn 3$ be a reduced curve
defined by the maximal minors of a $(4\times 5)$- matrix $M$ with linear
entries, let $\pi:V\to \Pn 3$ be the blow up of $\Pn 3$ along the ideal of
$C$, and let $\phi M: V\to \Pn 4$ be the morphism defined by the maximal minors
of $M$. Then any scheme $Z$ of length 2 on $V$ which
is mapped by $\phi M$ to a point is mapped by $\pi$  into a 4-secant line to
$C$.\par
$Proof.$ Any 4-secant line to $C$ is contracted by $\phi M$, so we
need to show the converse.  If $\pi (Z)$ does not have support on $C$, then it
follows from lemma $0.14$ that in the general plane $P$ through $\pi (Z)$,
there is a
4-secant line to $C$ through $\pi (Z)$, or conic through 8 points on $C$ or a
plane cubic through 10 points on $C$.  But no plane section of $C$ is contained
in a
cubic so the latter two is impossible.  If $\pi (Z)$ has length 2 and support
on
$C$, then the same argument works to show that $Z$ is embedded by $\phi M$
since $C$ has no 5-secant.$\slutt$ \medskip Let $S$ be the strict transform
of $\Sigma$ in $V$.  Then $S$ is smooth and the map $\phi M$ embeds $S$
into $\Pn 4$:  $S$ is normalized and at each singular point of the base curve
$\Sigma$ is of type $(0,1)$ or $(1,2)$, so by lemma 2.9, the strict
transform is
smooth over these points.  On the other hand, by lemma 2.11, any double
point of
the map comes from points on a 4-secant, and each 4-secant meets the quintic
surface $\Sigma$ once outside the base curve.  Clearly $S$ is a surface of
type B.

\bigskip
We proceed to the surfaces with $\chi =2$. From proposition $1.11$ we know that
$S$ has
three 6-secants in a plane, or one or no 6-secant.
 \proclaim Proposition 2.12.  If $\chi =2$ and $S$ has three
6-secants in a plane, then $S$ is birationally a K3-surface. It is the
image of a
complete intersection (2,3) in $\Pn 4$ with one quadratic singularity under
a map defined by quadrics through the singular point and tangent to the
surface at
three points which form a theta characteristic on a hyperplane section.\par
$Proof$.  Since $S$ has three 6-secants, $K^2=-3$ and $H\cdot K=6$ it
follows from the 6-secant formula that $S$ has three $(-1)$-conics, which
we denote by $E_1$, $E_2$ and $E_3$.  Thus $S$ is birationally K3 and
$K=E_1+E_2+E_3$.
By (1.11) the plane of the three 6-secants intersects the surface along a
quartic
curve $C$
with self-intersection $C^2=1$.  Let $|D|=|H-C|$ denote the residual
pencil.
\par \proclaim Lemma 2.13.  $C\cdot E_i=D\cdot E_i=1$, for $i=\overline
{1,3}$.\par
$Proof$. Now, $C\cdot E_1\leq 2$ since $E_1$ is a conic and cannot lie in
the plane of $C$.   Assume that $C\cdot E_1=2$.  Then $C\cup E_1$ spans a
hyperplane,
and the pencil $|D|$ contains a member $D=D_0+E_1$.  We get that $D_0\cdot
E_1=1$, while the arithmetic genus of $D$ is 4, so the arithmetic genus of
$D_0$ is also 4.  But $D_0$ has degree 5, so this is impossible by lemma $0.3$.
Similarly $C\cdot E_i\le 1$ for $i=2,3$.  Now, by adjunction, $C\cdot K=C\cdot
(E_1+E_2+E_3)=3$, so the lemma follows.$\slutt$\medskip

Let $S_1$ be the minimal model of $S$ and let $C_1$ and $D_1$ be the image on
$S_1$ of $C$ and $D$ respectively under the blowing down map.  Since $C\cdot
D=3$ it follows
from the lemma that $C_1\cdot D_1=6$ and that $C^2_1=4$ and $D^2_1=6$.  $D$ is
canonically embedded in $\Pn 3$, it is a complete intersection (2,3), so
the linear system $|D_1|$ defines a map $\phi {D_1}: S_1\to\Pn 4$, whose image
is a complete intersection $(2,3)$.  The image of $C_1$ is a curve of degree 6
and genus at least 3, so the image of $C_1$ is a hyperplane section.  Thus
$D_1-C_1$ is effective, it is a $(-2)$-curve, which we denote by $A$.  The
image
$\phi {D_1}(A)$ is a quadratic singularity.  Furthermore, since $C+A$ and $D$
are
linearly equivalent and the general $D$ does not meet $A$,  the
restrictions $\kn D{2D}$ and $\kn D{C+D}$ are isomorphic. But $\kn D{C+D}$ is
the
canonical bundle on $D$, so by adjunction $\kn D{C+D}\cong \kn
D{D+E_1+E_2+E_3}$.
Thus $\kn DD\cong \kn D{E_1+E_2+E_3}$ and $\omega_D\cong \kn
D{2(E_1+E_2+E_3)}$.
Therefore $E_1, E_2, E_3 $ are blown down to the points $p_1, p_2, p_3$ on
$S_1$,  which form a theta characteristic on the general $D_1$ through
these points.$\slutt$ \medskip

{\bf 2.14. Construction of C.}  For existence, let $V_3$ be a cubic
hypersurface in $\Pn 4$ with exactly one isolated quadratic singularity at a
point $x$,
and let $V_2$ be a smooth quadric hypersurface which contains the double point
of
$V_3$ such that the complete intersection $S_0=V_2\cap V_3$ has a quadratic
singularity at $x$ and is smooth elsewhere.   Let $$\pi_0 :X\to \Pn 4$$ be
the blowing-up of $\Pn 4$ in the point $x$, and let $S_1$ be the strict
transform of $S_0$ on $X$ with a curve $A$ lying over the point $x$.  Then
$S_1$ is
smooth, and the curve $A$ is an irreducible rational curve with
self-intersection
$A^2=-2$. Next, let $D_1$ be a general hyperplane section of $S_0$, and
let $\Pi$ be a plane in $\Pn 4$ which is tangent to $D_1$ in three distinct
points $p_1,p_2,p_3$. More specifically we require that the intersections
$V_2\cap\Pi$
and $V_3\cap \Pi$ are, respectively, an irreducible plane conic and a plane
cubic curve which both go through and have a common tangent at the points
$p_1,p_2,p_3$. Denote the preimage of the points $p_1,p_2,p_3$ on $S_1$
also by $p_1,p_2,p_3$, and blow them up with a map $\pi_1 :S\to S_1$ to get a
smooth surface $S$ with exceptional divisors $E_1,E_2,E_3$.
\par Let $D_1$ and $A$ also
denote the pullback (total transform) on $S$ of the hyperplane divisor on $S_0$
and of
the $(-2)$-curve over the node $x$ on $S_0$, respectively.
 Consider the linear system of curves $$|H|=|2D_1-A-\sum_{i=1}^3 2E_i|$$ on
$S$.

\proclaim Proposition 2.15.  The above data $V_2, V_3$ and $\Pi$ has been
chosen such that the linear system of curves $|H|$ on $S$ is very ample and
embeds
$S$ as a surface of degree 10 in $\Pn 4$.\par
$Proof$. The proof amounts to exploiting the decomposition of the divisor
$H$ given by $H\equiv C+D$, where
$$|D|=|D_1-\sum_{i=1}^3 E_i|$$ and  $$|C|=|D-A|.$$   Note that $|D|$ is a
pencil of curves on $S$, while $|C|$ contains exactly one curve, call it $C$,
namely the strict transform of the hyperplane section of $S_0$ which contains
the
points $p_1,p_2,p_3$ and $x$.

\proclaim Lemma 2.16.  The linear system $|H|$ restricts
to the canonical linear series on both the curve $C$ and any smooth
curve in $|D|$.  Furthermore the two restriction maps
$$\Coh 0SH\to \Coh 0DH$$ and
$$\Coh 0SH\to \Coh 0CH$$ are both surjective.\par
$Proof.$ For any smooth curve $D\in |D|$ the divisor
$(p_1+p_2+p_3)=(\sum_{i=1}^3 E_i) _{|D}$ is a theta-characteristic on $D$,
while $|D_1|$ is the adjoint linear system to $D$, so $$|K_D|=|{D_1}
_{|D}|=|(\sum_{i=1}^3 2E_i) _{|D}|.$$  Therefore $$
|H_{|D}|=|(2D_1-A-\sum_{i=1}^3 2E_i) _{|D}|=|(D_1-A)_{|D}|=|K_D|,$$ where
the latter equality holds since $A$ does not meet $D$ at all.\par  Next,
consider the cohomology of the exact sequences  $$\es SCSHDH$$ and  $$\es
SDSHCH$$
of sheaves on $S$.  As noted above, $\coh 0SC=1$, so by duality and
Riemann-Roch,
$$\coh 1SC =\coh 2SC =0.$$ Therefore, in the first sequence,  $$\coh 1DH =\coh
1D{K_D}=1$$ implies that $\coh 1SH =1$ and, by Riemann-Roch again, that
$\coh 0SH =5$. In the second sequence  $\coh 1SD =\coh 2SD =0.$ Thus $\coh 1SH
=1$
implies that $\coh 1CH =1$ and $\coh 0CH =3$.  Since $p_a(C)=3$ and $H\cdot
C=4$,
the linear series $|H|_C|$ is the canonical linear series on $C$.  Note
also that the restriction maps are both surjective.$\slutt$\smallskip

For the very ampleness of $|H|$, it suffices, by lemma 0.15, to check that
the linear series $|H_{|D}|$ is very ample for every curve $D\in |D|$.\medskip

 \proclaim Lemma 2.17.  $|H|$ restricts to a very ample linear series on
the curve $C+A$ on $S$.\par
$Proof.$  First check that $C$ is not
hyperelliptic.  As noted above, $C$ is the strict transform on $S$ of the
hyperplane section, call it $L_0$, of $S_0$ which contains the points
$p_1,p_2,p_3$ and $x$.  Since $L_0$ has a double point at $x$, we see that
the morphism defined by the canonical linear series on $C$ is the
projection of $L_0$ from the point $x$ into a plane.  Thus $C$ is hyperelliptic
if and only if this map is 2 to 1, or geometrically, any line in $\Pn 4$
through $x$ which meets $L_0$ in a point away from $x$ will meet $L_0$ in two
(possibly infinitely close) points away from $x$.  But by Bezout and our choice
of
$V_2$ and $V_3$,  this means that any such line is contained in the surface
$S_0$, which is absurd. Therefore $|H_{|C}|$ is very ample.\par
For $A$, note that $A\cdot H=2$, so $|H_{|A}|$ is very ample.  If
$\Coh 0SH\to \Coh 0AH$
is not surjective, or if $|H|$ does not separate points $p\in C\setminus A$
and $q\in A\setminus C$, then since $C\cdot A=2$, the curve $A$ is mapped
into the plane of $C$, which is absurd since by lemma 2.16 the image of
any curve in $|D|$ spans a $\Pn 3$.$\slutt$\par

\proclaim Lemma 2.18.  $|H|$ restricts to a very ample linear series on
the curves in $|D|$ which have one of the $E_i$'s as a component.\par

$Proof.$ If $D=G_1+E_1$ say, then $G_1\cdot E_1=2$, and $G_1$ is a plane
quartic curve, so if $|H|$ does not separate points $p\in G_1\setminus
E_1$ and $q\in E_1\setminus G_1$, then $|H_{|(G_1+E_1)}|$ must have
dimension 2, impossible.  For $G_1$ we get $$|(C+D)
_{|{G_1}}|=|(2D_1-\sum_{i=1}^3 2E_i)
_{|{G_1}}|=|(D_1-E_1)_{|{G_1}}|=|K_{G_1}|.$$
So $|H_{|{G_1}}|$ is very ample
unless $G_1$ is hyperelliptic.  But the hyperplane cutting $\pi (G_1)$ cuts
the quadric $V_2$ transversely at the point $\pi (E_1)$ since $\Pi\cap V_2$ is
smooth, so the projection from $\pi (E_1)$ is an isomorphism.
$|H_{|E_1}|$ is very ample since it has degree 2.  So it remains to show that
the
restrictions of sections are surjective in both cases.  For a general $D$ we
know that
the restriction is surjective, so $G_1+E_1$ spans a $\Pn 3$ and one of the
restrictions is surjective.  But if $G_1$ or $E_1$ is mapped to a line,
then $G_1+E_1$ lies in a plane since $E_1\cdot G_1=2$,
contradiction.$\slutt$\medskip

 For the rest of the curves $D$,  $|H_{|D}|$ is very ample
since the map $\pi :S\to S_0$ restricts to an isomorphism, $D_1$
is the pullback of the hyperplane divisor on $S_0$, and $$|(C+D)
_{|D}|=|2D _{|D}|=|(2D_1-\sum_{i=1}^3 2E_i) _{|D}|=|{D_1} _{|D}|.$$
This concludes the proof of proposition 2.15.$\slutt$\bigskip

\proclaim Proposition 2.19.  If $\chi =2$ and $S$ has one 6-secant,
then $S$ is birationally K3. It  has two $(-1)$-lines $E_1$ and $E_2$, and a
$(-1)$-curve $E_0$ of degree 4, and  two linear systems of curves
$|A_0|$ and $|B_0|$ of degree 6 and genus 4.  $S$ is embedded into $\Pn 4$
by the linear system of curves in $|A_0+B_0-4E_0-E_1-E_2|$ such that the
product
of the pencils of curves $|A_0-2E_0|$ and  $|B_0-2E_0|$ defines a map onto a
quadric surface, through which $E_0$ acquires one double point at a point
which coincides with the image of $E_1$ and $E_2$.\par

 $Proof.$ By the 6-secant formula $S$ has two $(-1)$-lines $E_1$ and
$E_2$, and since $K^2=-3$ and $H\cdot K=6$ the surface has a
$(-1)$-rational quartic curve $E_0$.  In fact the minimal model has a canonical
divisor
$K_0$ with $K_0^2\geq 0$ and $H\cdot K_0=0$ or $H\cdot K_0\geq 3$.  Only
$H\cdot
K_0=0$ fits in our case, so $S$ is birationally K3 and $K$ consists of the
three $(-1)$-curves.

\proclaim Lemma 2.20.  Any plane quartic curve $C$ on $S$ with
self-intersection $C^2=0$ will satisfy:  $$C\cdot E_0=2,\quad C\cdot
E_1=C\cdot E_2=1.$$\par
$Proof.$  Let $D=H-C$. Then $D$ has degree $H\cdot D=6$, arithmetic
genus 3 and self-intersection $D^2=2$.  Assume that $D\cdot E_i=e_i$ for
$i=\overline {0,2}$, then $D_0=D+\sum_{i=0}^2e_iE_i$ has self-intersection
$D_0^2=2+\sum_{i=0}^2e_i^2$.  Furthermore $D_0$ has arithmetic genus
$p_a(D_0)=3+\sum_{i=0}^2 {e_i\choose 2}$, and since $D_0\cdot E_i=0$,
 adjunction
gives $2p_a(D_0)-2=D_0^2$.  This leads to the relation $$e_0+e_1+e_2=2.$$
 Note
that $e_0>0$ since $H\cdot E_0=(C+D)\cdot E_0=4$ and $E_0$ spans all of $\Pn
4$; if
not there would be a curve $H-E_0$ on $S$ of degree 6 and arithmetic genus
5 which would necessarily have a plane quintic curve as a component,
impossible by  corollary 1.2.\par
Furthermore $e_i\geq 0$, for $i=1,2$, since $D$ moves in a pencil
with only isolated base points by corollary 1.2.  Therefore the only
possible
cases are $e_0=1$ and $\{e_1,e_2\}=\{0,1\}$, or $e_0=2$ and $e_1=e_2=0$.
The second case leads to the conclusion of the lemma.  In the first case assume
that $e_1=1$, the other case being analogous. Then there is a curve $D_2$ in
the
pencil $|D|$ which contains $E_2$ as a component.
 Thus $D_2=A+E_2$ where $A$ is a curve of degree 5 and genus 3.  This curve
again must decompose into a plane quartic curve $A_0$ and a line $L$ cutting
it.
So $A_0\cdot L=1$.  By the index theorem $A_0^2\leq 1$, and by lemma 1.3
$A_0^2\not= 1$, so $A_0^2\leq 0$.  Since $A^2=1$ this means that $L^2\geq -1$.
Thus $L^2 =-1$ and $L\cdot A=0$. But this holds only if $L=E_1$, which
leads to $E_1\cdot E_2=E_1\cdot (E_2+A)=E_1\cdot D=e_1=1$,
impossible.$\slutt$\medskip

\proclaim Lemma 2.21.  The hyperplane section of $S$ which contains the
exceptional lines decomposes into two plane quartic curves $A$ and $B$ with
$A\cdot B=2$ and $A^2=B^2=0$.  The 6-secant to $S$ is the line of
intersection between the planes of $A$ and $B$.\par
$Proof.$   First assume that the exceptional lines $E_1$ and
$E_2$ meet the 6-secant line, and let $C=H-E_1-E_2$, where $H$ is the
hyperplane
section containing $E_1$ and $E_2$.  By Riemann-Roch, $C$ lies on at least
two cubic surfaces, and the 6-secant line of $S$ intersect $C$ in a scheme of
length at least 4.  So if the $C$ lies on an irreducible cubic, it is linked to
the 6-secant line in the intersection of two cubics. Then the pencil of planes
through this line defines a pencil of degree 4 on $C$. The sum of this
pencil and a hyperplane section of $C$ form a canonical divisor on $C$.  But,
on
the $K3$ surface, $H+E_0$ restricts to $C$ to form a canonical divisor, which
means that $E_0$ intersect $C$ in a scheme of length 4 in a plane, impossible.
Therefore any cubic containing $C$ is reducible, and they all have a
component in common.\par
If this fixed component is a plane, then $C$ decomposes into a
plane quartic $A$ and an elliptic quartic curve $B$.  The arithmetic genus
of $C$ is 7, so $A\cdot B=4$.  Since  the 6-secant line is at least a 4-secant
to $C$ it must lie in the plane of $A$ and the pencil $|H-A|$ must have two
base points in the plane.  Thus $(H-A)^2=2$ and $A^2=0$.  But, by lemma 2.20,
this means that $A\cdot (B+E_1+E_2)=A\cdot (H-A)=6$, a contradiction.\par
If this fixed component is an irreducible quadric, then $C$ decomposes into a
line
$L$ and a curve $A$ on the quadric of arithmetic genus 6, with $L\cdot A=2$.
Again, the 6-secant line is a 4-secant to $C$ so it is a line on the quadric,
which is a 4-secant to $A$.  This is possible only if the intersection of
$L+E_1+E_ 2$
with the quadric has a scheme of length 2 residual to the intersection with
$A$,
i.e., that $A\cdot (E_1+E_2)=2$.  But this means that $L\cdot (E_1+E_2)=2$, so
$(H-L)\cdot L=4$ and $(L+E_1+E_2)^2=-1$, which is impossible.\par
So the fixed component must be two planes.  One of the planes must
intersect the surface along a quartic curve $A$ and the 6-secant must lie in
this
plane, so as above $A^2=0$ and $A\cdot E_i=1$.  Thus $(H-A-E_1-E_2)\cdot A=2$,
which
means that the curve $B=H-A-E_1-E_2$  must have arithmetic genus 3, so it
is also a plane quartic curve with self-intersection $B^2=0$, and the line
of intersection of the two planes is the 6-secant to the surface.\par
Now, assume that one of the exceptional lines, say $E_1$, does not
meet the 6-secant.    The 6-secant and $E_1$ span a hyperplane which
cuts the surface in $E_1$ and a curve $C_1$ of degree 9 and
arithmetic genus 8. Since any plane section of  $C_1$ through the 6-secant
lie on 3 cubic curves, and $\cohi 1{C_1}2=1$, $C_1$ must lie on a pencil of
cubics. This is possible only if the cubics have a fixed component.  The fixed
component cannot be a plane since $S$ does not contain any plane quintic curve.
 So the
fixed component is a quadric, and $C_1$ decomposes into a line $E$ and a curve
$C_0$, such that $C_0$ is on a quadric and $E\cdot C_0\leq 2$.  This means that
$E$ is a $(-1)$-line, i.e., $E=E_2$, and $C_0$ has arithmetic genus 7.  Thus
$C_0$
cannot lie on an irreducible quadric, and it is the union of two plane quartic
curves
$A$ and $B$ which meet in two points. Since $A\cdot E_1=A\cdot E_2=1$, it
follows that $B^2=0$, and similarly $A^2=0$.$\slutt$\medskip

    By (2.20) we have $A\cdot E_0=B\cdot E_0=2$.  Consider the two pencils
of curves $|A+E_1+E_2|$ and $|B+E_1+E_2|$, they have two base points each, so
the product of the pencils define a rational map of $S$ onto a quadric, this
map has degree 4 outside the base points, and the curve $E_0$ is mapped to an
elliptic quartic curve, so it acquires a double point.  Proposition 2.19 now
follows
from\par

 \proclaim Lemma 2.22. The inverse image of the double
point is two points on $E_0$ and the two $(-1)$-curves $E_1$ and
$E_2$.\par

$Proof.$  Since  $A+E_1+E_2=H-B$ and
$B+E_1+E_2=H-A$, the base points of the pencils $|A+E_1+E_2|$ and
$|B+E_1+E_2|$ lie in the planes of $B$ and $A$ respectively, in fact since $S$
has only
one 6-secant, the base points all lie on this 6-secant.  We denote the base
points of $|A+E_1+E_2|$ by $a_1$ and $a_2$ and the base points of $|B+E_1+E_2|$
by
$b_1$ and $b_2$.  There are two more points on the 6-secant which lie on the
surface, they form $A\cap B$ and we denote them by $c_1$ and $c_2$.\par

The curves $A$ and $B$ are canonically embedded on $S$ into $\Pn 4$. So
$a_1+a_2+c_1+c_2$ is a
canonical divisor on $A$.  But by adjunction the canonical line bundle on
$A$ is $\kn A{(A+E_1+E_2)+E_0}$, so $\kn A{A+E_1+E_2}=\kn A{a_1+a_2}$ means
that
$A\cap E_0=c_1+c_2$. Similarly $B\cap E_0=c_1+c_2$.  Now there is a net of
curves
in $|A+B+2E_1+2E_2|$ through the points $a_1, a_2, b_1, b_2, c_1, c_2$, and in
fact all these curves contain $E_1$ and $E_2$ so the lemma follows.$\slutt$\par
\medskip
For existence cf.(4.8).

 \proclaim Proposition 2.23.  If $\chi =2$ and $S$ has no 6-secants, then $S$
is
a smooth elliptic surface,
the canonical curve has four components: a plane cubic curve $C$ and three
$(-1)$-lines, $E_1$, $E_2$ and $E_3$. The linear system
$$|D|=|H+E_1+E_2+E_3-C|$$
defines a birational morphism of $S$ into $\Pn 3$ whose image is a surface
of degree 7 with a quadruple point at a point $q$, a double curve of degree 9
and arithmetic genus 9 lying on the cone with vertex at $q$ over a plane cubic
curve isomorphic to $C$, with 6 branches through $q$ along the 6 lines of
intersection of 4 planes through $q$.  Furthermore the surface $S$ lies
on a rational quartic with a double plane, the plane of the curve $C$.\par
$Proof.$  Since $S$ has
no 6-secant, it has three $(-1)$-lines and the canonical curve $K$ has the
components  $$K\equiv C+\sum_{i=1}^3E_i,$$ where $C$ is an elliptic curve of
degree 3, i.e., a plane cubic curve, and the $E_i$ are the $(-1)$-lines.
Thus $S$ is birationally an elliptic surface.\par   Consider the pencil
$$|D_0|=|H-C|$$ and the linear system $$|H_1|=|H+\sum_{i=1}^3E_i|$$
of curves on $S$.

\proclaim Lemma 2.24.  The
general member of $|D_0|$ is trigonal, and $|H_1|$  defines a morphism
which is the composition of the blowing down of the $(-1)$-lines and an
embedding
into $\Pn 6$. \par
$Proof.$ First note that the pencil $|D_0|$ has no fixed component,
since a fixed component of $|D_0|$  would be contained in the plane of $C$;
now $D_0\cdot C =3$, so this fixed component would be a line, call it $L$.
Thus $L\cdot C=3$ and $L\cdot 2C=6=H\cdot 2C$.  But $|2C|$ is a pencil of
elliptic curves whose general member is irreducible, so $L$ would be a 6-secant
to
such a curve, which is absurd. Thus the general member $D_0\in |D_0|$ is
irreducible.  It has degree $H\cdot D_0=7$ and arithmetic genus $p_a(D_0)=6$,
so it is linked to a line in the complete intersection of a quadric and a
quartic surface. One of the rulings on the quadric surface will sweep out a
$g^1_3$
on $D_0$, i.e., $D_0$ is trigonal.\par
 Furthermore $|H_1|=|D_0+C+\sum_{i=1}^3E_i|=|D_0+K|$ is the adjoint linear
system to $D_0$ on $S$. From the global sections of the exact sequence
$$\es SKS{H_1}{D_0}{H_1}$$ of sheaves of $S$, it follows, since $S$ is regular,
that the restriction map $\Coh 0S{H_1}\longrightarrow \Coh 0{D_0}{H_1}$ is
surjective and that $\coh 0S{H_1}=7$.  Thus $|H_1|$ is very ample on the
general
curve in $|D_0|$ and blows down the $E_i$'s. Since $|H|$ is very ample
the lemma follows.$\slutt$\par

\proclaim Corollary 2.25.  The embedding of $S$ into $\Pn 4$ is the
projection of a smooth surface in $\Pn 6$ from the linear span of
3 collinear points $x_1, x_2, x_3$.\par

$Proof.$  Denote the image of $S$ in $\Pn 6$ by $S_1$.  Let  $x_i=\phi
{H_1}(E_i)$, for $i=1,$ 2, 3.  Since $H\equiv H_1-\sum_{i=1}^3E_i$, the
points $x_i$ must be collinear. Moreover, they form a member of the
$g^1_3$ on the general $D_0$.$\slutt$\medskip

Denote by $L$ the line spanned by the points $x_1, x_2$
and $x_3$. Denote by $C_1$ the image of $C$ on $S_1$. The curve
$C_1$ is a plane cubic curve on $S_1$ whose plane we  denote by $\Pi$.
 Consider the linear system $$|D|=|H_1-C_1|$$ of curves on $S_1$. It has
projective
dimension dim$|D|=3$. Denote the subpencil of $|D|$  which corresponds to
the pencil $|D_0|$ on $S_1$ by ${\cal P}$; i.e., the pencil of curves in $|D|$
which meet the points $x_i$.  Although all of the curves in ${\cal P}$ are
trigonal, this is not necessarily the case for all the curves in $|D|$.

\proclaim Lemma 2.26.  There is a net of trigonal curves in $|D|$.\par

$Proof.$ To see this, consider the base locus
$Z_{\cal P}\subset S_1$ of the pencil ${\cal P}$.  We may write $Z_{\cal
P}=Z+Z_L$, where
$Z_L$ is of length three and is contained in $L$, while $Z$ is of length
four and has support outside $L$.  For a general member $D_p\in {\cal P}$, the
scheme $Z$ is a divisor which by duality on $D_p$ spans a $\Pn 3$ together with
the plane $\Pi$.  In fact $$Z+({\Pi}\cap {D_p})\equiv K_{D_p}-Z_L$$ as divisors
on $D_p$.  We denote this $\Pn 3 \subset \Pn 6$ by $V_{\cal P}$.  Now there is
a net
of curves $D$ in $|D|$ which contain $Z$, and, by duality again, all of these
curves are trigonal.$\slutt$\medskip

Denote this net by ${\cal N}$. Every member $D$ of ${\cal N}$ is a trigonal
curve canonically
embedded by $|H_1|$,  therefore it lies on a rational scroll denoted by $S_D$,
whose
ruling restricts to the trigonal linear series on $D$.  For the general member
$D_p\in {\cal P}$ the line $L$ must be a member of this
ruling.  $L$ cannot meet any other member of the ruling, otherwise the
projection of $D_p$ from the line $L$ into $\Pn 3$ would be three to one,
therefore the scroll $S_D$ must be smooth.\par\medskip

Let $V_0$ be the intersection of the quadric hypersurfaces containing $S_1$.
We proceed now to study $V$, the irreducible component of $V_0$
which contains $S_1$, and the projection of $V$ to $\Pn 3$
from the plane $\Pi$. Note that all rational normal scrolls in the net
$\{S_D| D\in {\cal N}\}$  are contained in $V_0$ since their rulings
are $3$-secants to $S_1$, and that $\Pi$ is also contained in
$V_0$.\par\medskip

First, let us consider the possibilities for $V_0$.  Since $\kn {H_1}{2H_1}$ is
non-special on any smooth curve $H_1$, it follows from Riemann-Roch and the
cohomology of the exact sequence $$\es {S_1}{H_1}{S_1}{2H_1}{H_1}{2H_1}$$
of sheaves on $S_1$ that $\cohi 0{S_1}2 =3$.\par

\proclaim Lemma 2.27. dim$V_0<4$.\par
$Proof$.  If $V_0$ is four-dimensional, then it is of degree three and
irreducible,
since it lies on three linearly independent quadrics.
Codimension two varieties of degree three are ruled by a pencil of linear
spaces of
codimension three. In our case $V_0$ is ruled by a pencil of $\Pn 3$s.
This pencil must clearly restrict to the ruling of the scrolls $S_D$.
Thus $L$ is contained in one of the $\Pn 3$s.
On the other hand, these $\Pn 3$s must sweep out a pencil of curves on
$S_1$.  Projecting from $L$, it follows that the member of this pencil which
belongs to the $\Pn 3$ of $L$, is mapped onto a line. Thus the curves of the
pencil must all be rational, which is absurd for an elliptic surface of
Kodaira dimension 1.$\slutt$\par\medskip

Thus $V_0$ and $V$ are three-dimensional and are contained in the complete
intersection of three quadric hypersurfaces. The scrolls $S_D$ are parts of
hyperplane sections of $V$:  The net ${\cal N}$ of divisors on $S_1$ is
the restriction to $S_1$ of the net of
hyperplanes which contains the linear space $V_{\cal P}$.  The restriction of
the
same net to $V$ has the scrolls $S_D$ as members of the moving part, since
the general such member must be irreducible. The fixed part of this net
is $T=V\cap V_{\cal P}$.

\proclaim Lemma 2.28.  $T$ is a cubic surface.\par
$Proof.$  Let $D_p$ be a general irreducible member of ${\cal N}$, and let
$S_D$ be the corresponding scroll.  First note that $S_D\cap V_{\cal P}$ is a
twisted cubic
curve; in fact if $E_0$ is a section with $E_0^2=0$ on $S_D$ and $F$ is a
member
of the ruling, then $D_p$ is of type $3E_0+4F$ on $S_D$, the adjoint linear
system is $|E_0+2F|$, so $|E_0+F|$ restricts to the $g^3_7$ residual to the
trigonal linear series defined by the ruling $|F|$.  Now $D_p\cap V_{\cal
P}=Z+({\Pi}\cap
{D_p})\equiv K_{D_p}-Z_L,$ is a member of this $g^3_7$, so it lies on a
curve in $|E_0+F|$ which is a twisted cubic curve.\par

\proclaim Lemma 2.29.  For a general member of ${\cal N}$ the corresponding
twisted cubic curve is irreducible.\par
$Proof.$  Assume that every twisted cubic is reducible.  Then the twisted
cubics
would have a line as a component which is a trisecant to  the corresponding
$D_p$.
Consider the plane of $C$ in $\Pn 4$.  The curves $|D_0|$ are of type $(3,4)$
on
a quadric and they have 4 base points in the plane.  Every trisecant to $D_p$
which is contained in $V_{\cal P}$ is a trisecant to the corresponding curve in
$|D_0|$ on $S$ and would lie in this plane and in the corresponding quadric.
Thus
all the quadrics meet the plane in two lines, and one of these lines is fixed
for the pencil of quadrics.  It would contain at least 3 base points so it
would be
a 6-secant to $S$, impossible.$\slutt$\par

 \proclaim Corollary 2.30.  $D_p$ and $S_D$  intersect $\Pi $ only on
$C_1$, and the scheme $Z$ spans $V_{\cal P}$.\par

$Proof$ $of$ $2.28$ $continued.$ The surface swept out by the twisted cubics
contains $C_1$, so it has a
component different from $\Pi $ which has degree at least 3.
Now, if $V_{\cal P}$ is not contained in $V_0$, then $V\cap V_{\cal P}$ is at
most a quadric surface.
It would have the plane $\Pi $ as a component, impossible by lemma 2.29.
Thus $V_{\cal P}$ is a component of $V_0$ which by linkage intersects the other
components in a cubic surface.  But $V$ intersects $V_{\cal P}$ in a surface of
degree at
least 3, so lemma 2.28 follows.$\slutt$.

 \bigskip

By (2.30), the plane $\Pi$ intersects each scroll in the scheme $C_1\cap S_D$
of
length 3, hence the linear system $|D|$ has no base points.
 The morphism $$\phi {D}:S_1\to \Pn 3$$ defined by $|D|$ is birational
since $D^2=7$ and the image spans $\Pn 3$.  The image $\Sigma =\phi
{D}(S_1)$  is a surface of degree seven.  The image of a
curve $D_p$ of the net ${\cal N}$ must be a plane curve of degree seven,
and therefore the image of the scroll $S_D$ by the projection from the plane
${\Pi}$ is a
plane in $\Pn 3$.  The base locus $Z$ of ${\cal N}$ is mapped onto a point
$$q=\phi {D}(Z).$$ Therefore $\phi {D}(D_p)$ acquires a quadruple point at
$q$.  Additionally, $\phi {D}(D_p)$ acquires three
double points from the members of the ruling of $S_D$ which meet $C_1$.
That there are no other singularities can be checked by the genus formula of a
plane curve.\bigskip

The morphism $\phi D$ extends to a rational map $$proj_{\Pi}:V \rto \Pn 3$$
which is generically finite.  Since $V$ is contained in the
complete intersection of three quadrics, it is birational. In fact,
if $P$ is a $\Pn 3\subset\Pn 6$ which contains the plane $\Pi$, then the
three quadrics will restrict to $P$ as the union of $\Pi$ and three other
planes.
 If the intersection of the other planes is finite, then it is one point.
 The family of scrolls $\{S_D\}$ is mapped by $proj_{\Pi}$ onto the net of
planes through
the point $q$, while the curves in ${\cal N}$ are mapped to plane septic curves
in these planes, with a quadruple point at $q$ and three extra double
points.\par\medskip

The image $C_0=\phi {D}(C_1)$ of the curve
$C_1$ is a plane cubic curve; it lies on a cubic cone with vertex at $q$,
which we denote by $S_3$. This cone is the image in $\Pn 3$ of the
exceptional divisor coming from the blowing up of $V$ along $C_1$.\smallskip

The inverse rational map
$$\rho:\Pn 3\rto V\subset \Pn 6,$$
restricts to a plane $proj_{\Pi}(S_D)$ as the map defined by a linear
system of plane quartic curves with a triple point at $q$ and three simple
points
outside $q$ as assigned base points, since the curve $D_p$ is canonically
embedded in $S_D$.  Therefore $\rho$ is defined by a linear system $|d_0|$
of quartic surfaces with a triple point, as assigned basepoint, at $q$, and
with
an assigned base curve, denoted by $C_B$, which meets a general plane through
$q$ in three points outside $q$.  \par
Let $$\pi:U\to U_0\to \Pn 3$$ be the composition of blowing up first $q$ to get
$U_0$ and then the strict transform of $C_B$ on $U_0$ to get $U$.
Let $E_{q,0}$ be the exceptional divisor on $U_0$, let $E_q$
be its strict transform on $U$, and let $E_B$ be the exceptional divisor
over $C_B$.  Then $\rho$ extends to a morphism  $$\rho_U:U\to V\subset
\Pn 6$$ which is defined by the linear system
$$|d|=|\pi^{\ast}d_0-3E_q-E_B|$$
of divisors on $U$.  The canonical divisor on $U$ is
$$K_U\equiv-\pi^{\ast}d_0 +2E_q+E_B.$$\medskip

The image by the map $\rho_U$ of the strict transform $E_q$ of $E_{q,0}\cong
\Pn 2$
is the cubic surface $T=V\cap V_{\cal P}$.  Now, $E_q$ is the
projective plane $E_{q,0}$ blown up in the points $q_i$, $i=\overline{1,n}$ of
the intersection with the
strict transform of $C_B$ on $U_0$.  The restriction to $E_q$ of the linear
system $|d|$ is the linear
system of plane cubic curves  on $E_{q,0}$ with assigned base points at
these points of intersection.
Therefore $T$ is a Del Pezzo surface, $n=6$, and $C_B$ has six branches at
$q$.  This means that $C_B$ is a curve of degree nine in $\Pn 3$.\par
Now, $\rho_U(U)=V\subset \Pn 6$ has degree 7.  On the other hand, the
restriction of $|d|$ to the strict transform of a general plane in $\Pn 3$
is the linear system of quartic curves with 9 assigned base points, which is a
linear system of degree 7 and projective dimension 5, so its image in $V$
must be a hyperplane section of  $V$. But this means that the base curve $C_B$
must be contained in $S_3$.\par\medskip
 Summing up:  The surface $\Sigma =\phi {D}(S_1)$ must have a quadruple
point at $q$ and must have $C_B$ as a double curve.  If $\Sigma_U$ is the
strict
transform of $\Sigma $ on $U$, then $\Sigma_U$ must meet $E_q$ in a curve
which is the strict transform of a quartic curve with double points at the
points
$q_i$, $i=\overline{1,6}$, on $E_{q,0}$.  Thus the points $q_i$ lie on a conic,
or
are the points of intersection of four lines in the plane $E_q$, where no three
lines meet in a point.  Let us exclude the
first of these cases.  In fact the strict transform  on $U$ of each plane
through $q$ intersects $E_q$ in the pullback of a line whose image in $V$ is a
twisted cubic curve.  If the 6 points lie on a conic, the image of the line in
$V$
would be a plane cubic curve and $S_1$ would be singular,
contradiction.\par
Thus we may consider the map of $\Sigma$ into $\Pn 4$ as defined by the
quartics containing the curve $C_B$ and a line through $q$.  Together, the
two curves form a curve $C_g$ of degree 10 and genus 11, which is linked
in the complete intersection of two quartics to a non-hyperelliptic
curve of genus 3 and degree 6. Therefore $C_g$ is defined by the maximal minors
of a
$4\times 5$-matrix with linear entries.  The 5 minors define the map to $\Pn 4$
and
the image of $\Pn 3$ is a determinantal quartic.  Moreover, the exceptional
surface over
$q$ is mapped two to one onto the plane spanned by $C$, so the quartic
is singular along this plane.$\slutt$ \medskip

{\bf 2.31. Construction of E}.
It is clear from the above description how to construct such
a surface.  Start with a smooth plane cubic through the 6 points of
intersection of four lines. On the cone over this curve with vertex $q$, look
at curves linear equivalent to the plane curve plus the lines through the 6
special
points.  Such a curve is linked on the cone via a
quartic with triple point at $q$ to
three lines through the vertex, which span $\Pn 3$, so it is clear from Bertini
that one can find a curve $C_B$ in the above linear system which is smooth
outside the vertex $q$.  Let $U$ denote the blow up, first of $\Pn 3$ at $q$
and then of the strict transform of $C_B$.  Let  $S_{3,U}$ and $E_q$ denote
the strict transforms on $W$ of the cubic cone and of the
exceptional divisor over $q$, respectively.  If $h$ denotes
the pullback of a plane, and $E_B$ denotes
the exceptional divisor over $C_B$, then $S_{3,U}$ belongs to the linear system
$|3h-E_B-3E_q|$ of divisors on $U$.  The linear system
$$|d|=|h+S_{3,U}|=|4h-E_B-3E_q|$$ of divisors on $U$ defines the map
$\varphi_{|d|}:U\to\Pn 6$. Consider the linear system of divisors
$$|\Sigma_U|=|7h-2E_B-4E_q|=|d+S_{3,U}+2E_q|$$ on $W$.
  Observe that $${(\Sigma_U)}_{|E_q}\equiv L_1+L_2+L_3+L_4.$$ Furthermore, the
linear system $|\Sigma_U|$ has base points only on $S_{3,U}$. In fact
restricting to
this scroll one can easily show that there are no base points outside $E_q$.
Therefore, by Bertini, one can choose a smooth member
$\widetilde{S_1}\in |\Sigma_U|$, making sure that
$$(\widetilde{S_1}\cap S_{3,U})\cap (\widetilde{S_1}\cap E_q)=\emptyset.$$
The canonical divisor of $\widetilde{S_1}$ is given by adjunction:
$$K_U\equiv -4h+2E_q+E_B,$$ so $$K_{\widetilde{S_1}}\equiv
(K_U+\widetilde{S_1})_{|\widetilde{S_1}}
\equiv (S_{3,U}+E_q)_{|\widetilde{S_1}}\equiv
(\widetilde{S_1}\cap S_{3,U})+L_1+L_2+L_3+L_4.$$
 Therefore $L_1+\dots+L_4$ is part of the canonical divisor
on $\widetilde{S_1}$. Now if $1\le i\le 4$, then
$$K_{\widetilde{S_1}}\cdot L_i=L^2_i$$
since $L_i$ does not meet any of the other components of
$K_{\widetilde{S_1}}$.  Thus, the curves $L_i$, $i=\overline{1,4}$,
are $(-1)$-curves on $\widetilde{S_1}$, which are blown down on
$S_1=\varphi_{|d|}(\widetilde{S_1})$.\par\medskip
Pick a general line through $q$. Together with $C_B$ it forms a
curve of degree $10$, genus $11$, which is defined by the $4$-minors
of a $4\times 5$ matrix with linear entries. The line is mapped to a
line $L$  in $\Pn 6$, which is a trisecant of $S_1$. The composition
of $\varphi_{|d|}$ and the projection from $L$ into $\Pn 4$ is given
by the quartic minors above, and maps $S_1$ onto a surface $S$.
As in the case of the rational surface (cf. 2.10) the map defined by these
quartics  would get double points only from the 4-secants to the base curve.
Therefore it is easy to check that $S$ is embedded in $\Pn 4$.\par

\proclaim Proposition 2.32.  If $S$ is a smooth
surface of degree 10 in $\Pn 4$ with $\pi =9$ and $\chi =3$, then
$K^2=3$, $p_g=2$, $q=0$, and $S$ has exactly one $(-2)$-curve $A$ such that $S$
is embedded in $\Pn 4$ by the linear system $|2K-A|$.\par

$Proof.$ We first study the pencil of canonical curves $|K|$.  Now, $H\cdot
K=6$
and $p_a(K)=K^2+1=4$, so a general integral curve in the pencil must be
canonically embedded by $|H|$.  Thus, if $K$ is a general canonical curve
and we consider the cohomology of the exact sequence $$\es
S{K-H}S{2K-H}K{2K-H}$$
of sheaves on $S$, then we must have that $\coh 0K{2K-H}=1$.  On the other
hand $\coh 1S{K-H}=\coh 1SH=0$ by the Severi's theorem and Riemann-Roch, so we
get
that $\coh 0S{2K-H}=1$.  Let $A$ be the curve of $|2K-H|$.  Then $H\cdot A=2$
and $K\cdot A=0$ and $A^2=-2$, so $A$ is a (possibly reducible) $(-2)$-curve of
degree two on $S$.$\slutt$\medskip

{\bf 2.33. Construction of F}.
For existence consider the image $S_0$ of $S$ by the bicanonical map.
Then $S_0\subset \Pn 5$.  If $x,y$ are linearly independent sections of
$\kn S{K}$, then there is a quadratic relation between the sections
$x^2,xy,y^2$ of $\kn S{2K}$, therefore $S_0$ lies on a quadric $Q$ of rank
3.
 Let $\rho : X\to Q$ be the desingularization map.
$X=\openP(E)$ is a $\Pn 3$-bundle over
a $\Pn 1$.  On $X$ there is a divisor corresponding to the section of the
bundle
$E\otimes \kn {\Pn 1}{-2}$, call it $B$, and let $F$ denote a fiber of the
natural projection $\pi :X\to \Pn 1$.  $B$ is contracted to the vertex
plane of
$Q$ by $\rho$, in fact the map $\rho$ is defined by the linear system of
divisors
$|B+2F|$. We will construct a surface $S_1$ on $X$ which is mapped onto
$S_0$ in
$Q$ by $\rho$.  In fact $S_1$ is linked in the complete intersection of two
divisors $D_1\in |2B+6F|$ and $D_2\in |3B+6F|$  to three quadric surfaces
 $Q_0+Q_1+Q_2$, which are fibers of
the projection $\pi :D_1\to \Pn 1$.  Now, $S_0$ has one quadratic
singularity;
this is imposed by a quadratic singularity on $D_2$, which is cut
transversally by $D_1$, and which does not meet $B$ or any of the quadric
surfaces $Q_i$.   A Bertini argument assures the smoothness of $S_0$
outside
the singularity, and that $S_0$ has a quadratic singularity. The projection
from the singular point to $\Pn 4$ induces an embedding as soon the fiber
$Q_p$
 of  $D_1$ which meets the singularity is smooth.\medskip
\proclaim Remark 2.33.  An upshot of this is that the projection of
$D_2\subset Q$ to $\Pn 4$ is a quartic hypersurface with a plane, which
contains
$S$.\par
In fact the pencil of surfaces $D_2\cap F$ are mapped to cubic surfaces
which are residual to a plane in hyperplane sections of the quartic.

\proclaim Proposition 2.34.  If $S$ is a smooth surface of degree 10 with
$\pi=10$ in $\Pn 4$, then $S$ is a regular elliptic surface  with two
$(-1)$-lines and $p_g=2$; it is the projection from a point $p$ in $\Pn 5$
of a surface $S_0$ which is linked $(2,3,3)$ to a Del Pezzo surface of
degree 6, such that both $S_0$ and the Del Pezzo has an improper node at
$p$, or $S$ is a minimal regular surface of general type with
three $(-2)$-curves $A_1,A_2$ and $A_3$ embedded as
conics and $p_g=3$, such that $S$ is embedded by the linear system
$|2K-A_1-A_2-A_3|$ in $\Pn 4$.\par
$Proof.$  By proposition 1.16 and remark
1.17, it follows that if $\chi =3$, then $S$ has no 6-secant and two
$(-1)$-lines, and if $\chi =4$ then $S$ has one 6-secant and no $(-1)$-lines.
For the case $\chi =3$, consider the map defined by the linear system
$|H_0|=|H+E_1+E_2|$, where $E_1$ and $E_2$ are the $(-1)$-lines on $S$.
Clearly, this linear system defines a birational morphism of $S$ into $\Pn 5$
which sends the exceptional lines to a point $p$.  We denote by $S_0$ the
image. Then $S_0$ has an improper node at $p$ and the embedding of $S$
into $\Pn 4$ is the projection from this point.\par

\proclaim Lemma 2.35. $S_0$ lies on one quadric and on three
independent cubic hypersurfaces.\par
$Proof.$  Let $C_H$
be a general hyperplane section through $p$, and consider the cohomology of
the exact sequence  $$\esi {S_0}1{S_0}2{C_H}2.$$ Since $S_0$ is linearly normal
in $\Pn 5$ the surface $S_0$ lies on a quadric as soon as $C_H$ does.  But
$\kn {C_H}2$ is non-special, so $C_H$ lies on a quadric by Riemann-Roch. On the
other hand, if $S_0$ lies on more then one quadric then the projection  of a
complete intersection of two quadrics through $p$ yields a quadric or a cubic
hypersurface containing $S$, impossible. A similar argument with the above
exact
sequence, using also a general section not through $p$, gives the result for
cubics.$\slutt$\par\medskip

\proclaim Lemma 2.36.  $S_0$ is contained in a complete intersection
$(2,3,3)$.\par
$Proof.$  Let $V$ be the intersection of all the cubics containing $S_0$.
If $V$ is a threefold, then $V$ has degree at most 4,
since it is contained in four cubics which are independent of the quadric
$Q$. But then the projection of $V$ from the point $p$ has degree at most 3 and
contains $S$, impossible.$\slutt$\par\medskip

 Let $T$ be the surface which is linked to $S_0$ in a general
pencil of cubics containing $S$ on the quadric $Q$.  Since $p_g(S_0)=2$ and
the general canonical curve is an elliptic normal curve of degree 6, the
surface $T$ lies on a pencil of irreducible quadric hypersurfaces. Since $T$
has
degree 6 it is linked $(2,2,2)$ to a surface $U$ of degree 2.  Using a formula
analogous to (0.11) on $Q$ it follows that $T$ has sectional genus 1 and
$U$ has sectional genus -1. Locally at $p$, both $T$ and $U$ must have an
improper node, so $U$ must be two planes which meet at $p$. From the minimal
free resolution of $U$, we get a minimal free resolution of $S_0$.
In particular we see that $S_0$ is cut out by the quadric and the 4 cubic
hypersurfaces. Therefore, by Bertini, the general $T$ is smooth outside
the node. It is the projection of a Del Pezzo in $\Pn 6$ from
a point on a secant line, such that the image
$T$ in $\Pn 5$ acquires an improper node.\par\medskip
  For the case $\chi =4$,  consider a general member $C_K\in |K|$. Since
$H\cdot C_K=8$ and $p_a(C_K)=K^2+1=5$, the linear series $|H_{|C_K}|$ is not
canonical on
$C_K$, only if $C_K$ spans only a hyperplane in $\Pn 4$.  But $p_g=\coh 0SK\geq
3$, so the residual curve $D\equiv H-C_K$ has degree $H\cdot D=2$ and is
contained in a line, which is absurd.  Therefore $C_K$ is canonically
embedded
in $\Pn 4$ by $|H_{|C_K}|$.   Consider now the cohomology associated to the
exact sequence $$\es S{K-H}S{2K-H}{C_K}{2K-H}$$ of sheaves on $S$.
By Severi's theorem and the Riemann-Roch theorem
$\coh 1S{K-H}=\coh 1SH=0$ and $\coh 0S{K-H}=0$,
and by the above $\kn {C_K}{2K}\cong \kn {C_K}H$, so $\coh
0{C_K}{2K-H}=1$. Therefore $\coh 0S{2K-H}=1$.  Let $A$ be the curve in
$|2K-H|$.  Since $K\cdot A=0$, $A$ must be the union of $(-2)$-curves on $S$.
Now $A^2=-6$ and $p_a(A)=-2$, so $A$ must be the union of three numerically
disjoint $(-2)$-curves which we denote by $A_1, A_2$ and $A_3$.  Thus
$A_i\cdot
A=A_i\cdot (2K-H)=-2$, for $i=\overline{1,3}$.  But $K\cdot A_i=0$, so we get
that
$H\cdot A_i=2$, which means that the $(-2)$-curves $A_1,A_2$ and $A_3$ are
embedded
as conics in $\Pn 4$.\par\medskip
Therefore $S$ is a regular surface embedded by the
linear system $$ |2K-A_1-A_2-A_3|,$$ in $\Pn 4$. This concludes the
proof of Proposition 2.34.$\slutt$\par\bigskip

For existence, the elliptic surface
can be reconstructed via the linkage in $\Pn 5$, but we will give a different
proof of existence using linkage in $\Pn 4$ (cf. 4.15).\par\medskip

For the surface of general type the easiest proof of existence is by using
linkage
in $\Pn 4$, (cf. 4.20). Alternatively, consider the image $S_1$ of $S$ by the
bicanonical map. Since the sections of  $\kn S{2K}$ contain the square of the
sections of $\kn SK$, the surface $S_1$ must lie on a four-dimensional cone
$X$ over a Veronese surface in $\Pn 7$.  Let  $X_0$ be the natural
desingularization
of $X$.  Then  $$X_0={\bf P}(E)={\bf P}(\kn {\Pn 2}{2}\oplus{\cal O}_{\Pn 2}
\oplus{\cal O}_{\Pn 2}).$$
Let $\rho : X_0\to X$ be the desingularization map.
On $X_0$ there is a divisor corresponding to
the section of the bundle $E\otimes \kn {\Pn 2}{-2}$, call it $B_0$, and
let $F$ denote the pullback
of a line by the natural projection $\pi :X_0\to \Pn 2$.  $B_0$ is
contracted to the vertex line
of $X$ by $\rho$, and in fact the map $\rho$ is defined by the linear system of
divisors $|B|=|B_0+2F|$.
$X_0$ has a canonical divisor $K_0\equiv -3B-F$, and we compute for the
intersection numbers $B^4=4$
and $B^3\cdot F=2$ and $B^2\cdot F^2=1$.  For two general members $G_1$ and
$G_2$ of $|2B|$, the complete intersection $G_1\cap G_2$ is a smooth surface
$\Sigma$.  By adjunction, $B-F$ restricts to the canonical divisor $K_{\Sigma}$
on $\Sigma$.  Thus $\Sigma$ is a regular surface with $\chi ({\Sigma}) =4$
and $K^2_{\Sigma}=4$.  Since $B^3\cdot B_0=0$, $\Sigma $ is embedded by $\rho$
in $\Pn 7$ as the complete intersection of $X$ and two quadric hypersurfaces
$Q_1$ and $Q_2$. The surface  $S_1$ is such a complete intersection with three
linearly independent quadratic singularities.  The map to $\Pn 4$ is the
projection from the plane spanned by these nodes.
\par\bigskip\bigskip

{\bf 3\quad Syzygies}\bigskip

 To describe the minimal free resolution of the ideal sheaf of a surface $S$ in
$\Pn 4={\bf P}(V)$
we use Beilinson's spectral sequence [Bei] whose $E_1$-terms are
$$E^{pq}_1=H^q({\bf P}(V),F(p))\otimes \Omega^{-p}(-p)$$ and which converges to
$$E^m=\cases{F,&for $m$=0;\cr 0,& otherwise;\cr}$$ where $F=\id Sm$ for a
suitable m. In the
sequel, the proper twist $m$ will always be 4.\par
Furthermore for the differentials we have the following:
\proclaim Lemma [Bei] 3.1. There are canonical isomorphisms
$Hom(\Omega^{i}(i),\Omega^{j}(j))\cong
\Lambda^{i-j} V$ defined by contraction and the composition of morphisms
coincides with
multiplication in $\Lambda V$.\par
\proclaim Lemma [D] 3.2. Let $A$ be a $s\times t$ matrix with entries in $V$
defining the
morphism $$A: \dsum\limits_t \Omega^{i}(i) \; \longrightarrow \;
\dsum\limits_{s}\Omega^{i-1}(i-1).$$
If $A$ is pointwise surjective then for any nontrivial linear combination of
the rows of $A$
$(a_1,...,a_{t})$ we have that ${dim}_k$ $span_{V}(a_1,\dots,a_{t})\geq i +
1$.\par

The irreducibility claims in the table at the beginning of the paper will be
consequences of the
following statement:\par

\proclaim Lemma  3.3. Let ${\cal E}$ and ${\cal F}$ be two vector bundles on
$\Pn 4$ of ranks $r$ and
respectively $r+1$ and let $\varphi_1$ and $\varphi_2$ be two injective
morphisms between ${\cal E}$ and
${\cal F}$. Assume also that both determinantal loci $S_i=\lbrace\; p\in\Pn 4\;
\mid\; {\rm rk}\,
{\varphi}_i (p)< r\;\rbrace$, $i=1,2$ have the expected codimension two. Then
$S_1$ and $S_2$
lie in the same irreducible component of the Hilbert scheme.\par
$Proof.$  Let $X = {\rm Spec} (k\lbrack t\rbrack)$ and let $\varphi =
(1-t)\varphi_1 + t\varphi_2$. One deduces
easily the existence of an open subset of $X$ containing $0$ and $1$ and
parameterizing a flat deformation
from $S_1$ to $S_2$. See also [K] or [BoM].$\slutt$\par\bigskip

The results of this section will recover some of the constructions via the
Eagon-Northcott complex
performed in [DES].\par

\proclaim Proposition 3.4. The ideal sheaf of a smooth rational surface $S$ of
degree 10 in $\Pn 4$,
with $\pi=9$, embedded by
$$H\equiv8\pi^{\ast}\li -\sum_{i=1}^{12}2E_i-\sum_{j=13}^{18}E_j$$
has a minimal free resolution:
$$
\vbox{%
\halign{&\hfil$\,#\,$\hfil\cr
&&&&2{\cal O}(-4)\cr
&&&&\oplus\cr
0&\leftarrow&{\cal I}_S&\leftarrow&5{\cal O}(-5)&&9{\cal O}(-6)&&3{\cal
O}(-7)\cr
&&&&\oplus&\vbox  to
10pt{\vskip-4pt\hbox{\nwarrow}\vss}&\oplus&\longleftarrow&\oplus\cr
&&&&{\cal O}(-6)&&3{\cal O}(-7)&&3{\cal O}(-8)&\vbox  to
10pt{\vskip-4pt\hbox{\nwarrow}\vss}&{\cal O}(-9)&
\leftarrow 0.\cr
}}$$\par
$Proof.$ The surface has one 6-secant and $\deh 4$ is the plane corresponding
to all hyperplanes through
the unique 6-secant. Lemmas 1.4, 1.15 and proposition 1.11 give also the
cohomology diagram:
\test \ \ 2 \ \ 2 1 2
$$\beili55v \ \ 2 \ \ 2 1 2 {$h^i({\cal I}_S (p))$}$$
where empty boxes are zeroes. Beilinson's spectral sequence provides then a
presentation of the ideal sheaf:
$$ 0\longrightarrow \; 2\Omega^{3}_{\Pn 4}(3)\; \longrightarrow \; {\cal F}
\oplus 2{\cal O} \; \longrightarrow \; {\cal I}_S(4) \; \longrightarrow 0 \eqno
\hbox{($\ast$)}$$
with ${\cal F}$ the kernel of the differential $d_1$:
$$ 0\longrightarrow \; {\cal F}\; \longrightarrow \; 2\Omega^{1}_{\Pn 4}(1) \;
\mapright {d_1} \; {\cal O}\; \longrightarrow 0$$ See [DES] for a construction
of the surface $S$
using this presentation. From lemma 3.2 it follows that the morphism $d_1$ is
defined by two linearly
independent elements of $V$ so an easy computation yields a minimal free
resolution for the vector
bundle ${\cal F}$:
$$
\vbox{%
\halign{&\hfil$\,#\,$\hfil\cr
&&&&15{\cal O}(-1)&&11{\cal O}(-2)&&3{\cal O}(-3)\cr
0&\leftarrow&{\cal
F}&\longleftarrow&\oplus&\longleftarrow&\oplus&\longleftarrow&\oplus\cr
&&&&{\cal O}(-2)&&3{\cal O}(-3)&&3{\cal O}(-4)&\vbox  to
10pt{\vskip-4pt\hbox{\nwarrow}\vss}&
{\cal O}(-5)&\leftarrow 0\cr}}
\eqno \hbox{($\ast \ast$)}$$
The claim of the proposition follows now from ($\ast$) and ($\ast \ast$).
The three linear forms in the last syzygy above, or equivalently
part of the last syzygy in the minimal resolution of ${\cal I}_S$,
annihilate the module $I_S/(I_S)_{\le 5}$, where $I_S$ denotes the homogeneous
ideal of $S$. Hence this module is supported
on a line, which thus is the unique $6$-secant line to $S$, the dual
of the variety $\deh 4$.$\slutt$\par

\proclaim Proposition 3.5. The ideal sheaf of a smooth rational surface of
degree 10 in $\Pn 4$,
with $\pi=9$, embedded by $$H\equiv 9\pi^{\ast}\li
-\sum_{i=1}^43E_i-\sum_{j=5}^{11}2E_j-\sum_{k=12}^{18}E_k$$
has a minimal free resolution:\par
$$
\vbox{%
\halign{&\hfil$\,#\,$\hfil\cr
&&&&{\cal O}(-4)\cr
0&\leftarrow&{\cal I}_S&\leftarrow&\oplus\cr
&&&&10{\cal O}(-5)&\vbox to 10pt{\vskip-4pt\hbox{\nwarrow}\vss}&18{\cal
O}(-6)&\leftarrow&10
{\cal O}(-7)&\leftarrow&2{\cal O}(-8)\leftarrow 0.\cr
}}$$\par
$Proof.$ From proposition 1.11, remark 1.12 and the proposition 2.1 it follows
that $S$ has no 6-secant line, and that $\cohi 1S4\le 1$.  If $\cohi 1S4=1$,
then
$\deh 4$ is a plane, namely the set of hyperplanes through a line $L$, lying
this time on $S$ since the surface has no $6$-secants. But then the argument
in proposition 3.4 would show that $I_S/(I_S)_{\le 5}$ is supported on $L$,
which in particular means that $L$ is not on $S$, a contradiction.
Together with lemmas 1.4 and 1.15 we obtain therefore the
following cohomology table:
\test \ \ 2 \ \ 2 \ 1
$$\beili55v \ \ 2 \ \ 2 \ 1 {$h^i({\cal I}_S (p))$}$$ Beilinson's
spectral sequence recovers a presentation of the ideal sheaf as in [DES]:
$$ 0\longrightarrow \; 2\Omega^{3}_{\Pn 4}(3)\; \longrightarrow \;
2\Omega^{1}_{\Pn 4}(1)
\oplus {\cal O} \; \longrightarrow \; {\cal I}_S(4) \; \longrightarrow 0 $$
and consequently the announced resolution. $\slutt$\par

\proclaim Proposition 3.6. The ideal sheaf of a smooth elliptic surface $S$ of
degree 10 in $\Pn 4$,
with $\pi=9$ has a minimal free resolution:
$$
\vbox{%
\halign{&\hfil$\,#\,$\hfil\cr
&&&&{\cal O}(-4)\cr
0&\leftarrow&{\cal I}_S&\leftarrow&\oplus\cr
&&&&9{\cal O}(-5)&\vbox  to 10pt{\vskip-4pt\hbox{\nwarrow}\vss}&14{\cal
O}(-6)&&5{\cal O}(-7)\cr
&&&&&&\oplus&\longleftarrow&\oplus\cr
&&&&&&{\cal O}(-7)&&2{\cal O}(-8)&\vbox  to
10pt{\vskip-4pt\hbox{\nwarrow}\vss}&{\cal O}(-9)&
\leftarrow 0\cr
}}.$$\par
$Proof.$ It follows from lemmas 1.4, 1.15 and proposition 1.11 that
$S$ has no 6-secants
and that $\deh 4$ is a line, namely the set of all
hyperplanes containing the plane $P$ of the elliptic cubic curve which is part
of the canonical
divisor. Moreover, lemma 1.4, 1.14 and proposition 1.15 give the following
cohomology table:
\test 1 \ 1 \ 1 3 1 1
$$\beili55v 1 \ 1 \ 1 3 1 1 {$h^i({\cal I}_S (p))$}$$ whence Beilinson's
spectral sequence provides
the presentation:
$$ 0\longrightarrow \; {\cal O}(-1)\oplus\Omega^{3}_{\Pn
4}(3)\oplus\Omega^{2}_{\Pn 4}(2)\;
\longrightarrow \; {\cal F} \oplus {\cal O} \;
\longrightarrow \; {\cal I}_S(4) \; \longrightarrow 0 \eqno \hbox{($\ast$)}$$
where ${\cal F}$ is the kernel of $d_1$:
$$ 0\longrightarrow \; {\cal F}\; \longrightarrow \; 3\Omega^{1}_{\Pn 4}(1) \;
\mapright {d_1} \; {\cal O}\; \longrightarrow 0.$$
Since $\deh 4$ is a line, $d_1$ is given by three linearly independent elements
of $V$ so we get the
minimal resolution of ${\cal F}$:
$$
\vbox{%
\halign{&\hfil$\,#\,$\hfil\cr
0&\leftarrow&{\cal F}&\leftarrow&25{\cal O}(-1)&&&20{\cal O}(-2)&&6{\cal
O}(-3)\cr
&&&&&&\vbox  to
10pt{\vskip-4pt\hbox{\nwarrow}\vss}&\oplus&\longleftarrow&\oplus\cr
&&&&&&&{\cal O}(-3)&&2{\cal O}(-4)&\vbox  to
10pt{\vskip-4pt\hbox{\nwarrow}\vss}&{\cal O}(-5)&
\leftarrow 0\cr
}}.\eqno \hbox{($\ast \ast$)}$$
where the two linear forms of the last syzygy above are defining the plane $P$.
{}From ($\ast$) and ($\ast \ast$)
we get easily the minimal resolution of the ideal sheaf of $S$.$\slutt$\par

\proclaim Proposition 3.7. The ideal sheaf of a smooth $K3$ surface $S$ of
degree 10 in $\Pn 4$,
with $\pi=9$ and one 6-secant has a minimal free resolution:
$$
\vbox{%
\halign{&\hfil$\,#\,$\hfil\cr
&&&&{\cal O}(-4)\cr
&&&&\oplus\cr
0&\leftarrow&{\cal I}_S&\leftarrow&9{\cal O}(-5)&&15{\cal O}(-6)&&7{\cal
O}(-7)&&{\cal O}(-8)\cr
&&&&\oplus&\vbox  to
10pt{\vskip-4pt\hbox{\nwarrow}\vss}&\oplus&\longleftarrow&\oplus
&\longleftarrow&\oplus&\leftarrow 0\cr
&&&&{\cal O}(-6)&&3{\cal O}(-7)&&3{\cal O}(-8)&&{\cal O}(-9)\cr
}}.$$\par
$Proof.$ We proceed as for the elliptic surface. In this case the cohomology
diagram is the same
but $\deh 4$ is now a plane, namely the set of all hyperplanes through the
6-secant line $L$.
Beilinson's spectral sequence gives then a presentation:
$$ 0\longrightarrow \; {\cal O}(-1)\oplus\Omega^{3}_{\Pn
4}(3)\oplus\Omega^{2}_{\Pn 4}(2)\;
\longrightarrow \; {\cal F} \oplus {\cal O} \;
\longrightarrow \; {\cal I}_S(4) \; \longrightarrow 0 \eqno \hbox{($\ast$)}$$
with ${\cal F}$  the kernel of $d_1$:
$$ 0\longrightarrow \; {\cal F}\; \longrightarrow \; 3\Omega^{1}_{\Pn 4}(1) \;
\mapright {d_1} \; {\cal O}\; \longrightarrow 0.$$
The morphism $d_1$ is given by a triple $(v_1, v_2, v_3)\in V^3$, but since
$\deh 4$ is a plane,
$v_1, v_2$ and $v_3$ are linearly dependent so we can assume that $v_3=0$ and
$v_1 , v_2$ are linearly independent.
Therefore ${\cal F} = {\cal F}^{\prime}\oplus\Omega^{1}_{\Pn 4}(1)$ where
 $${\cal F}^{\prime} =
{\rm ker}(2\Omega^{1}_{\Pn 4}(1)\maponto {\hbox{$\scriptstyle(v_1, v_2)$}}
{\cal O}).$$ Moreover, using
the minimal free resolution of ${\cal F}^{\prime}$ we computed in proposition
3.4, one gets for ${\cal F}$:
$$
\vbox{%
\halign{&\hfil$\,#\,$\hfil\cr
&&&&25{\cal O}(-1)&&21{\cal O}(-2)&&8{\cal O}(-3)&&{\cal O}(-4)\cr
0&\leftarrow&{\cal
F}&\longleftarrow&\oplus&\longleftarrow&\oplus&\longleftarrow&\oplus
&\longleftarrow&\oplus&\leftarrow 0\cr
&&&&{\cal O}(-2)&&3{\cal O}(-3)&&3{\cal O}(-4)&&{\cal O}(-5)\cr
}}$$
and hence from ($\ast$) the desired resolution of the ideal sheaf of
$S$.$\slutt$\par

\bigskip
As a consequence of the above two propositions we get that, for both elliptic
and $K3$ surfaces,
the $H^1({\cal I}_S(\ast))$ -module is a monogeneous artinian module with
Hilbert function $(1,3,1)$
over $R = k{\lbrack}x_0,\dots x_4{\rbrack}$, but the module structure is
different. We use the above
description to sketch a construction of such surfaces. We produce first a
vector bundle ${\cal G}$
of rank 6 with the property that $H^1({\cal G}(\ast)) = H^1({\cal I}_S(\ast +
4))$ and $H^i({\cal G}(\ast)) = 0$
for $i = 2,3$. Let then ${\cal F} = {\cal O}(-1)\oplus\Omega^{3}_{\Pn 4}(3)$.
In both cases
the degeneracy locus of a generic $\varphi\in{\rm Hom}({\cal F},{\cal G})$
will be a surface with the desired invariants. Also the Eagon-Northcott complex
of $\varphi$
$$ 0\longrightarrow \; {\cal F}\; \longrightarrow \; {\cal G}\;
\longrightarrow\; {\cal I}_S(4)\; \longrightarrow 0.$$
will lead to the above computed resolutions. In both cases we take
$${\cal G} = {\rm ker}(5{\cal O}\oplus2{\cal O}(1)\mapontoa {\psi} {\cal
O}(2)),$$
where the linear part of $\psi$ is given, without loss of generality, by $x_0$
and $x_1$ and the
quadratic by $q_1,\dots q_5\in k{\lbrack}x_2,x_3,x_4{\rbrack}$ quadrics in
three variables without
common zeroes. Therefore the choice of $\psi$ is equivalent to that of a
hyperplane section of the
Veronese surface in $\Pn 5$. But there are only two possible types: an
irreducible hyperplane section
(generic one) whose corresponding quadrics define a ${\cal G}$ leading to the
elliptic surface, or a
reducible hyperplane section - two conics with a common point - whose quadrics
define a ${\cal G}$
leading to the $K3$ surface. In the first case ${\rm dim}_k {\rm Hom}({\cal
F},{\cal G})=35$,
while in the second case it is $36$, whence we easily deduce that both families
of surfaces have
dimension $44$.\par
We remark also that a variation of the hyperplane section of the Veronese
surface gives a deformation
of elliptic surfaces of degree 10, $\pi=9$ to a scheme belonging to the
irreducible component of the
Hilbert scheme containing the $K3$ surfaces of the above type.\par

\proclaim Proposition 3.8. The ideal sheaf of a smooth $K3$ surface $S$ of
degree 10 in $\Pn 4$,
with $\pi=9$ and with three 6-secants has a minimal free resolution:
$$
\vbox{%
\halign{&\hfil$\,#\,$\hfil\cr
&&&&2{\cal O}(-4)\cr
&&&&\oplus\cr
0&\leftarrow&{\cal I}_S&\leftarrow&4{\cal O}(-5)&&7{\cal O}(-6)&&2{\cal
O}(-7)\cr
&&&&\oplus&\vbox  to
10pt{\vskip-4pt\hbox{\nwarrow}\vss}&\oplus&\longleftarrow&\oplus\cr
&&&&3{\cal O}(-6)&&8{\cal O}(-7)&&7{\cal O}(-8)&\vbox  to
10pt{\vskip-4pt\hbox{\nwarrow}\vss}&2{\cal O}(-9)&
\leftarrow 0\cr
}}.$$\par
$Proof.$ Lemmas 1.4, 1.15 and proposition 1.11 give the following cohomology
table:\par
\test 1 \ 1 \ 1 3 2 2
$$\beili55v 1 \ 1 \ 1 3 2 2 {$h^i({\cal I}_S (p))$}.$$
so Beilinson's spectral sequence gives then a presentation of ${\cal I}_S(4)$:
$$ 0\longrightarrow \; {\cal O}(-1)\oplus\Omega^{3}_{\Pn
4}(3)\oplus\Omega^{2}_{\Pn 4}(2)\;
\longrightarrow \; {\cal F} \oplus 2{\cal O} \;
\longrightarrow \; {\cal I}_S(4) \; \longrightarrow 0 \eqno \hbox{($\ast$)}$$
where ${\cal F}$ is the kernel of $d_1$:
$$ 0\longrightarrow \; {\cal F}\; \longrightarrow \; 3\Omega^{1}_{\Pn 4}(1) \;
\mapright {d_1} \; 2{\cal O}\; \longrightarrow 0.$$
Moreover, by lemma 1.6 and proposition 1.11, $\deh 4$ is a degenerated cubic
scroll: three planes with a common
line. By lemma 3.1, $d_1$ is defined by a matrix $ D = \left( \matrix{
v_1&v_2&v_3\cr w_1&w_2&w_3\cr}
\right)$ with entries in $V$ whose $2\times2$ minors, when we regard it as a
matrix of linear forms
in $\Pnd 4$, are cutting out the scroll $\deh 4$. But $\deh 4$ is degenerated
so making, if
necessary, a change of coordinates $D$ can be brought, by lemma 3.2, to the
form
$ D = \left( \matrix{ u&v&0\cr 0&v&w\cr}\right)$ with $u, v, w$ linearly
independent elements of $V$.
Hence we can compute the minimal free resolution of $\cal F$:
$$
\vbox{%
\halign{&\hfil$\,#\,$\hfil\cr
&&&&20{\cal O}(-1)&&13{\cal O}(-2)&&3{\cal O}(-3)\cr
0&\leftarrow&{\cal
F}&\longleftarrow&\oplus&\longleftarrow&\oplus&\longleftarrow&\oplus\cr
&&&&3{\cal O}(-2)&&8{\cal O}(-3)&&7{\cal O}(-4)&\vbox  to
10pt{\vskip-4pt\hbox{\nwarrow}\vss}&
2{\cal O}(-5)&\leftarrow 0\cr
}}.$$
The proposition follows now easily. We remark also that the plane spanned by
the 6-secant lines to
$S$ is ${\bf P}({span}_k (u, v, w))$ and, more precisely, that the three
6-secants to $S$ are in fact
${\bf P}(\overline{uv}) , {\bf P}(\overline{uw})$ and ${\bf
P}(\overline{vw})$.$\slutt$\par
\proclaim Remark 3.9. {\rm Conversely, to construct such a surface one takes
$${\cal G} = {\cal O}\oplus{\rm ker}(4{\cal O}\oplus2{\cal O}(1)\mapontoa
{\psi} {\cal O}(2))$$
with $\psi$ generic and ${\cal F} = {\cal O}(-1)\oplus\Omega^{3}_{\Pn 4}(3)$.
For a generic $\psi\in
{\rm Hom}({\cal F},{\cal G})$ the degeneracy locus $S = \{\,x\in \Pn 4 \mid{\rm
rk}\varphi(x) <
{\rm rk}{\cal F}\,\}$ is a $K3$ surface of the considered type. See also [DES]
for this construction}.\par

\proclaim Proposition 3.10. The ideal sheaf of a smooth general type surface
$S$ of degree 10 in
$\Pn 4$, with $\pi=9$ has a minimal free resolution:
$$
\vbox{%
\halign{&\hfil$\,#\,$\hfil\cr
&&&&{\cal O}(-4)\cr
&&&&\oplus\cr
0&\leftarrow&{\cal I}_S&\leftarrow&8{\cal O}(-5)&&13{\cal O}(-6)&&6{\cal
O}(-7)&&{\cal O}(-8)\cr
&&&&\oplus&\vbox  to
10pt{\vskip-4pt\hbox{\nwarrow}\vss}&\oplus&\longleftarrow&\oplus
&\longleftarrow&\oplus&\leftarrow 0\cr
&&&&3{\cal O}(-6)&&8{\cal O}(-7)&&7{\cal O}(-8)&&2{\cal O}(-9)\cr
}}.$$\par
$Proof.$ From lemma 1.4, 1.15 and proposition 1.11 we get the cohomology table:
\test 2 \ \ \ 2 4 2 1
$$\beili55v 2 \ \ \ 2 4 2 1 {$h^i({\cal I}_S (p))$}$$
and that $\deh 4$ is a degenerated cubic scroll: three planes having a common
line, corresponding to
the three 6-secants of $S$. As in the above proposition we obtain the
presentation:
$$ 0\longrightarrow \; 2{\cal O}(-1)\oplus2\Omega^{2}_{\Pn 4}(2)\;
\longrightarrow \; {\cal F} \oplus {\cal O} \oplus \Omega^{1}_{\Pn 4}(1) \;
\longrightarrow \; {\cal I}_S(4) \; \longrightarrow 0 $$
with ${\cal F}$ the kernel of the map given by the matrix $D$ in the proof of
proposition 3.8:
$$ 0\longrightarrow \; {\cal F}\; \longrightarrow \; 3\Omega^{1}_{\Pn 4}(1) \;
\mapright D \; 2{\cal O}\; \longrightarrow 0.$$
The claim of the proposition follows, and we get a similar description of the
6-secant lines.$\slutt$\par

\proclaim Remark 3.11. {\rm Conversely to construct such a surface we take
$${\cal G} = {\rm ker}
({\cal O}\oplus6{\cal O}(1)\mapontoa {\psi} 2{\cal O}(2))$$ for a special
morphism : the $2\times 6$
matrix which is the linear part of $\psi$ should drop rank in exactly three
points. Let ${\cal F} = 2{\cal O}
(-1)\oplus 2{\cal O}$. Then for a generic $\varphi\in {\rm Hom}({\cal F},{\cal
G})$ one obtains a surface of
general type of degree 10, with $\pi=9$. We remark also that the three 6-secant
lines to the surface are
exactly the lines joining the points where the linear part of $\psi$ is
dropping rank}.\par

\proclaim Proposition 3.12. The ideal sheaf of a smooth elliptic surface $S$ of
degree 10 in $\Pn 4$, with $\pi=10$ has a minimal free resolution:
$$
\vbox{%
\halign{&\hfil$\,#\,$\hfil\cr
&&&&3{\cal O}(-4)\cr
0&\leftarrow&{\cal I}_S&\longleftarrow&\oplus\cr
&&&&3{\cal O}(-5)&\vbox  to 10pt{\vskip-4pt\hbox{\nwarrow}\vss}&9{\cal O}(-6)&
\longleftarrow&5{\cal O}(-7)&\longleftarrow&{\cal O}(-8)\leftarrow 0.\cr
}}
$$\par
$Proof.$ In proposition 1.16 we determined the following cohomology table:
\test 2 \ 1 \ \ 1 \ 3
$$\beili55v 2 \ 1 \ \ 1 \ 3 {$h^i({\cal I}_S (p))$}$$
Therefore we obtain a presentation of the ideal sheaf as cokernel of:
$$ 0\longrightarrow \; 2{\cal O}(-1)\oplus\Omega^{3}_{\Pn 4}(3)\;
\longrightarrow \; 3{\cal O} \oplus \Omega^{1}_{\Pn 4}(1) \;
\longrightarrow \; {\cal I}_S(4) \; \longrightarrow 0 \eqno \hbox{($\ast$)}$$
hence the claimed resolution. See [DES] for a construction, using $(\ast)$, of
the
surface $S$ .$\slutt$\par

\proclaim Proposition 3.13. The ideal sheaf of a smooth, general type surface
$S$ of
degree 10 in $\Pn 4$, with $\pi=10$ has a minimal free resolution:
$$
\vbox{%
\halign{&\hfil$\,#\,$\hfil\cr
&&&&3{\cal O}(-4)&&{\cal O}(-5)\cr
&&&&\oplus&&\oplus\cr
0&\leftarrow&{\cal I}_S&\leftarrow&3{\cal O}(-5)&\longleftarrow&6{\cal
O}(-6)&&2{\cal O}(-7)\cr
&&&&\oplus&&\oplus&\vbox  to 10pt{\vskip-4pt\hbox{\nwarrow}\vss}&\oplus\cr
&&&&{\cal O}(-6)&&3{\cal O}(-7)&&3{\cal O}(-8)&\vbox  to
10pt{\vskip-4pt\hbox{\nwarrow}\vss}
&{\cal O}(-9)&\leftarrow 0\cr
}}.$$\par
$Proof.$ From proposition 1.16 it follows that $S$ has one 6-secant, that $\deh
4$
is formed by all hyperplanes containing the 6-secant and the following
cohomology
table:
\test 3 \ \ \ 1 2 1 3
$$\beili55v 3 \ \ \ 1 2 1 3 {$h^i({\cal I}_S (p))$}.$$
Beilinson's spectral sequence gives then a presentation of the ideal sheaf of
$S$ as cokernel of:
$$ 0\longrightarrow \; 3{\cal O}(-1)\oplus\Omega^{2}_{\Pn 4}(2) \;
\longrightarrow \; {\cal F} \oplus 3{\cal O} \;
\longrightarrow \; {\cal I}_S(4) \; \longrightarrow 0 \eqno \hbox{($\ast$)}$$
with ${\cal F}$, as in 3.4, the kernel of the differential $d_1$:
$$ 0\longrightarrow \; {\cal F}\; \longrightarrow \; 2\Omega^{1}_{\Pn 4}(1) \;
\mapright {d_1} \; {\cal O}\; \longrightarrow 0  \eqno \hbox{($\ast \ast$)}$$
Moreover the natural map $$H^1({\cal I}_S (2))\otimes H^0({\cal O}_{\Pn
4}(1))\longrightarrow H^1({\cal I}_S (3))$$
is surjective and $\deh 3 = \deh 4$. The claimed resolution follows now easily.
{}From the description
of $\deh 3$ we get also that the linear syzygy of the three quartic generators
of the ideal sheaf of
$S$ is formed by the equations of the unique 6-secant line. The same remark
applies for the linear
part of the last syzygy.$\slutt$\par\smallskip
Conversely, as above, one uses $(\ast)$ and $(\ast\ast)$ to
reconstruct the surface.\par\bigskip\bigskip

 {\bf 4\quad Linkage}\bigskip

 In this section we are going to describe minimal elements of the liaison
classes to
which the above studied surfaces belong and, therefore, to provide alternative
constructions for some of the surfaces. For general facts concerning even
liaison classes of codimension 2, locally Cohen-Macaulay subschemes of $\Pn n$
see [PS], [LR], [BoM] and [MDP]. We'll need also in the sequel the following
version of a lemma
from [LR]:\bigskip

\proclaim Lemma 4.1. Let $Z$ be a codimension two, locally Cohen-Macaulay
subscheme of $\Pn n$,
and define the speciality of $Z$ as $e(Z)$:=$\max {\{\,t\mid h^{n-2}({\cal
O}_{Z}(t)) \not= 0\,\}}$.\hfil\break
a) If $h^{0}({\cal I}_{Z}(e(Z) + n )) = 0$, then $Z$ is a minimal element in
its even liaison class.\hfil\break
b) If, moreover, $h^{0}({\cal I}_{Z}(e(Z) + n + 1)) = 0$, then $Z$ is the
unique minimal element
in its even liaison class.\par
$Proof.$ The lemma is stated in [LR] only for space curves but the proof works
also in the
general case.$\slutt$\par\bigskip

\proclaim Corollary 4.2. A smooth rational surface $S$ of degree 10 in $\Pn 4$,
with $\pi=9$, embedded by
$$H\equiv 9\pi^{\ast}\li
-\sum_{i=1}^43E_i-\sum_{j=5}^{11}2E_j-\sum_{k=12}^{18}E_k$$
is minimal in its even liaison class.\par\bigskip

\proclaim Proposition 4.3. A smooth rational surface $S$ of degree 10 in $\Pn
4$, with $\pi=9$,
embedded by $$H\equiv8\pi^{\ast}\li -\sum_{i=1}^{12}2E_i-\sum_{j=13}^{18}E_j$$
is minimal in its even liaison class and is linked (4,4) to a reducible surface
$Z=T_1\cup T_2$
where $T_2$ is a smooth cubic scroll and $T_1$ is a degenerated cubic scroll,
union of three
planes having in common the directrix line of $T_2$ and such that there are no
further
intersection points with $T_2$. Moreover $Z$ is the unique minimal scheme in
its even liaison
class.\par
$Proof.$ Associated with the surface $S$, we got in proposition 2.2 three plane
quartic curves $C_{12},
C_{34}$ and $C_{56}$ such that their corresponding planes $\Pi_{12},\Pi_{34}$,
resp. $\Pi_{56}$
have in common a line $L$, the unique 6-secant line to $S$. Moreover, $S$ meets
each of the planes
$\Pi_{ij}$ in two points on $L$ outside the curve $C_{ij}$ so we get that each
such plane has
two pencils of 5-secant lines to $S$. Now the family of 5-secants to $S$ is one
dimensional
since $S$ is contained in a pencil of irreducible quartics. Therefore, by  Le
Barz's formula,
we obtain that $$T_1=\Pi_{12}\cup\Pi_{34}\cup\Pi_{56}$$ is the union of the
5-secants to $S$.\par
$S$ can be linked (4,4) to a surface $Z$ of degree 6 with $\pi=1$. $T_1$ is
contained in any
quartic hypersurface containing $S$ so it will necessarily be a component of
$Z$. Consider now
the blowing up $p:\tilde S\to S$ of $S$ in the points $\{ q_1,\dots q_6 \} =
S\cap L$. We denote by $G_1,\dots,
G_6$ the exceptional divisors. The projection of $S$ along $L$ defines a
morphism of degree 4
$\varphi :\tilde S\to \Pn 2$ which is in fact the morphism induced by ${\tilde
D}\equiv p^{\ast}H - {\sum_{i=1}^6 G_i}$.
Recall now from the end of the proof of proposition 2.2 the existence of a
curve $L_1\in\mid{\pi}^{\ast}\li\mid$
passing through all the points $q_i$ , $i=\overline{1,6}$ and consider then
$\tilde{L_1} =
{(p\circ\pi)}^{\ast}\li - \sum_{i=1}^6G_i$. $\tilde{L_1}$ is mapped to a plane
conic by $\varphi$
since $\tilde D\cdot\tilde{L_1} = 2$ and obviously $\mid\tilde D -
\tilde{L_1}\mid = p^{\ast}\mid
H - {\pi}^{\ast}\li\mid = \emptyset$. Therefore there exists a curve $\tilde
C\equiv 2\tilde D -
\tilde{L_1}$ on $\tilde S$ and we denote by $C$ its image on $S$; ${\rm deg}\,
C = 12$ and
$p_a (C) = 12$.\par

\proclaim Lemma 4.4. $C$ lies on a smooth cubic scroll $T_2$ containing $L$ as
directrix.\par
$Proof.$ We have $h^2({\cal O}_S( 2H - C)) = h^0({\cal O}_S(K + C - 2H)) = 0$
so Riemann Roch on $C$ and
the cohomology of the exact sequence: $$\eso C 2 S 2 C 2$$ give that $h^0({\cal
I}_C(2))\ge 2$.
Since $h^0({\cal I}_C(1)) = 0$ it follows that also $h^0({\cal I}_{C \cap
H}(2))\ge h^0({\cal I}_C(2))\ge 2$ for a generic hyperplane $H$. Let $D =
p(\tilde D)$ for a general $\tilde D\in\mid\tilde D\mid$. Then $D$ meets $C$ in
six points on $L$ and in six points $x_1,\dots x_6$
on $S$ outside $L$.\par We work for a moment on $\tilde S$. $\varphi$ restricts
to a map of degree 3 on $\tilde C$, so we may group the $x_i$'s on $S$ into two
sets such that, say $x_1, x_2, x_3$
span a plane which contains $L$ and $x_4, x_5, x_6$ span another plane which
contains $L$.
Since $\varphi(\tilde C)$ is a conic we get that $x_1 + x_2 + x_3$ and $x_4 +
x_5 + x_6$ belong to the
same linear series on $C$. Now $L\cup \{ x_1,\dots x_6\}$ is contained in at
least two quadrics
so we get that either $\{ x_1, x_2, x_3\}$ or $\{ x_1, x_2, x_3\}$ is contained
in a line. But they
belong, as divisors on $C$, to the same linear series so they are both
contained in a line. Varying
$D$ we see that the lines $L_1$ and $L_2$ are describing the ruling of a scroll
$T_2$ containing the
curve $C$ and the line $L$. Moreover $T_2$ is contained in all quadrics
containing $C$ and is smooth
since $L_1$ and $L_2$ cannot meet. The scroll $T_2$ is rational and has a
hyperplane section divisor
$H_{T_2}\equiv 2\li - E$, viewed as $\Pn 2$ blown up in one point. The line $L$
equals the directrix
$E$, while $C$ meets $L$ in six points, so $C\cdot E = 6$. Since deg $C = 12$
we get on $T_2$:
$C\equiv 9\li - 6E.$ $\slutt$\par\bigskip

{}From the above lemma it follows that any conic $\li$ on $T_2$ cuts $C$ in 9
points and therefore any
quartic hypersurface that contains $S$ must contain also $T_2$. It follows that
$T_2$ is a component of $Z$
whence $$Z = T_1\cup T_2 = \Pi_{12}\cup\Pi_{34}\cup\Pi_{56}\cup T_2$$ as
claimed. For the rest of the
assertions one computes $e(S) = -1$, $e(Z) = -2$ and uses lemma
4.1.$\slutt$\par\bigskip

{\bf 4.5. Construction of surfaces A.} We start now with two cubic scrolls as
above, $T_1$
and $T_2$ such that $T_1$ is the union of three planes $\Pi_1$, $\Pi_2$ and
$\Pi_3$ through a
common line $L$, and $T_2$ is smooth having $L$ as directrix and meeting $T_1$
only along $L$.
The surface $Z=T_1\cup T_2$ has degree 6 and sectional genus $\pi=1$.\par

\proclaim Lemma 4.6. The scheme $Z$ is cut out by $11$ quartic hypersurfaces
and the generic
quartic in $I_Z$ is singular along $L$ containing it with multiplicity $2$.\par
$Proof.$ Consider the residual exact sequences
$$\esi {\Pi_3\cup T_2} 3 Z 4 {\Pi_1\cup\Pi_2\cup L'\cup L'',H_1} 4,$$
$$\esi {T_2} 2 {\Pi_3\cup T_2} 3 {\Pi_3\cup M'\cup M'',H_2} 3,$$
where $H_1$ is the hyperplane spanned by $\Pi_1\cup\Pi_2$, $L'$ and $L''$ are
rulings
of $T_2$ such that $L\cup L'\cup L''=H_1\cap T_2$, while $H_2$ is a general
hyperplane through
$\Pi_3$ and $M'$, $M''$ are other two rulings of $T_2$ defined by $H_2\cap T_2=
L\cup M'\cup M''$. The lemma follows since $\cohi i {T_2} 2=0$, for $i=1,2$,
$T_2$ is defined
by the minors of a $2\times 3$ matrix with linear entries, $\id {\Pi_3\cup
M'\cup M'',H_2} 3
\simeq\id {M'\cup M'',H_2} 2$ and $\id {\Pi_1\cup\Pi_2\cup L'\cup L'',H_1}
4\simeq
\id {L'\cup L'', H_1} 2$ are globally generated, while $\cohi 1 {M'\cup
M'',H_2} 2=0$.$\slutt$
\par\bigskip

As a consequence of the above lemma, $Z$ can be linked in the complete
intersection of two
general quartic hypersurfaces to an irreducible surface $S$ with the desired
invariants:
$\deg\, S=10$, $\pi=9$, $p_g=q=0$. It is easily seen that, for general choices
in the linkage,
$S$ is smooth outside the line $L$. To see the behavior at the intersection
points with $L$
we'll work out explicitly this linkage. \par\bigskip

Consider the blowing-up
$$ \sigma : \widetilde{\Pn 4}=\openP(2\sstruct{\Pn 2}\oplus\struct{\Pn 2} 1)
\longrightarrow\Pn 4$$
of $\Pn 4$ along the line $L$, with exceptional divisor $E=\openP(2\sstruct{\Pn
2})=
\openP(3\sstruct L)=\Pn 2\times L$. Let $B$ be a divisor of $E$ corresponding
to a section
of $3\struct L 1$ and $F$ corresponding to a fibre of the projection $\sigma :
\Pn 2\times L \to L$. If a hypersurface $V$ of degree $v$ contains the line $L$
with
multiplicity $m$, then its strict transform $\widetilde V$ will meet $E$ along
$\overline V$,
numerically of type
$$ \overline{V}\equiv mB + (v - 2m)F.$$
\par\bigskip

Let now $V_1$ and $V_2$ be two general quartic hypersurfaces through $Z$. By
lemma $4.6$,
they have multiplicity two along $L$, hence $\overline{V_i}\equiv 2B$,
$i=\overline{1,2}$.
On the other hand the strict transforms of the scroll $T_2$ and of a plane
$\Pi_i$ cut
$E$ along $\overline{T_2}\equiv B^2$, since the cubic scroll is linked to a
plane in the
complete intersection of two hyperquadric cones simple along $L$, and along
$\overline{\Pi_i}\equiv (B - F)(B - F)\equiv B^2 - 2BF$, respectively.
Therefore
$\overline{T_1}\equiv 3(B^2-2BF)$, and for a general $(4,4)$ linkage the strict
transform of $S$ on $\widetilde{\Pn 4}$ meets $E$ in a curve equivalent to
$(2B)(2B) - (B^2 + 3B^2 - 6BF)\equiv 6BF$. A Bertini argument shows now that
for a general
choice of the linkage, the surface $S$  is smooth, and thus it is a rational
surface as claimed.
Moreover, since a curve of type $BF$ is blown down on $S$, it follows that $L$
is a
$6$-secant line to $S$. We remark also that the quartics containing $S$ are
singular
along the $6$-secant $L$.\par\bigskip

\proclaim Proposition 4.7. A smooth $K3$ surface $S$ of degree 10 in $\Pn 4$,
with $\pi=9$ and three
6-secants is linked $(4, 4)$ to a reducible surface $Z = P_1\cup P_2\cup
P_3\cup T$, where $T$ is a
multiplicity three structure on a plane $P$ given, after a suitable change of
coordinates, by the
homogeneous ideal $$I_T = (\,{x}_0^2, {x_0}{x_1}, {x}_1^3, a{x}_1^2 +
{b_1}{b_2}{b_3}x_0\,)$$ with
$a, b_i \in k{\lbrack}x_2,x_3,x_4{\rbrack}, i=\overline{1,3}$, homogeneous
forms of degree 2, resp. 1,
without common factors when assuming $P = \{\,x_0 = x_1 = 0\,\}$, and where
$P_1$, $P_2$ and $P_3$
are planes which pairwise span all of $\Pn 4$ and cut the scheme $T$ along
multiplicity two structures
on the lines $L_i = \{\,x_0 = x_1 = b_i = 0\,\}$, for $i=\overline{1.3}$.
Moreover, the scheme $Z$ is
minimal in its even liaison class and is linked $(3,4)$ to a scheme $Y$ of
degree 6, with $\pi=1$,
which is minimal in the even liaison class of $S$.\par

$Proof.$ We shall use in the sequel notations from proposition 2.12: $C$ is the
plane quartic curve and
$P$ its plane, $D\equiv H - C$ the residual pencil, $\pi: S\to S_1$ the
blowing down of the three $(-1)$-conics $E_1$, $E_2$ and $E_3$ to the points
 $p_1$, $p_2$ and resp. $p_3$.  $D_1$ is the image of $D$  on
$S_1$, and corresponding ${\varphi}_{\mid{D_1}\mid}: S_1\to {\Pn 4}$ the
induced morphism whose image $S_0$
is a $(2,3)$ complete intersection with one quadratic singularity. By
proposition $2.12$, the points
$p_1$, $p_2$, $p_3$ form a theta-characteristic on the general $D_1$ through
these points. Since $D_1$
is canonically embedded, it follows that $p_1$, $p_2$, $p_3$ and the tangent
plane to $S_0$ at $p_i$
span only a hyperplane in $\Pn 4$, for all $i=\overline{1,3}$. Therefore we
obtain the existence of a
curve $G_i\in\mid D - E_i\mid$, for all $i=\overline{1,3}$ (cf. also $2.14$).
It has degree 4 and
arithmetic genus 3 so it is a plane curve. We denote by $P_i$ the plane spanned
by $G_i$.
Let $q_1, q_2, q_3$ be the base points of the pencil $D$ in the plane $P$ and
let $L_{ij}=\overline{q_i q_j}$
be the three 6-secant lines to $S$. Now $G_i\cdot C = 2$, $E_i\cap C =\{q_i\}$
and $D\cap C =
\{\,q_1,q_2,q_3\,\}$, hence $P_i$ cuts $P$ along the line $L_{jk}$, for all
triples $\{\,i,j,k\,\} =
\{\,1,2,3\,\}$. Since $S$ meets each of the planes $P_i$ in two points on
$L_{jk}$ outside the curve $G_i$,
for $\{\,i,j,k\,\} = \{\,1,2,3\,\}$, we get that each plane $P_i$ has two
pencils of 5-secants lines to $S$
and therefore is contained in any quartic hypersurface containing $S$.\par
Now the surface $S$ is linked $(4,4)$ to a surface $Z$ of degree 6, with
$\pi=1$ which must have as
components the planes $P_i$, $i=\overline{1,3}$. We consider further the scheme
$X = S\cup {P_1}\cup
{P_2}\cup {P_3}$. It has degree 13 and sectional genus 18 but it is not locally
Cohen-Macaulay, having exactly
six bad points $m_1, n_1, m_2, n_2, m_3, n_3$, where $L_{ij}\cap C =
\{\,q_i,q_j,m_k,n_k\,\}$ for all
triples $\{\,i,j,k\,\} = \{\,1,2,3\,\}$. Therefore $X$ is linked in the
complete intersection of the
above two quartic hypersurfaces to a scheme $T$ of degree 3, with $\pi=-2$. The
generic hyperplane
section of $T$ is a locally $CM$, connected curve of degree 3 and arithmetical
genus -2, hence it is
a triple structure on a line, hence $T$ also is a multiplicity three structure
on a plane. But, by
Bezout's theorem, the plane $P$ is contained in the complete intersection of
 the two quartic hypersurfaces,
whence $T$ must be a triple structure on $P$. Moreover, the liaison exact
sequences give that each
plane $P_i$ cuts the scheme $T$ along a doubling of the 6-secant line $L_{jk}$,
for all $\{\,i,j,k\,\} =
\{\,1,2,3\,\}$. Now a similar computation to [Ma] yields that, may be after a
suitable change of
coordinates, we may assume, when $P = \{\,x_0 = x_1 = 0\,\}$, that $T$ is given
by the homogeneous
ideal $I_T = (\,{x}_0^2, {x_0}{x_1}, {x}_1^3, a{x}_1^2 + b x_0\,)$ with $a,
b\in k{\lbrack}x_2,x_3,x_4{\rbrack}$
homogeneous forms of degree 2, resp. 3 such that ${\rm codim}_P(\,V(a)\cap
V(b)\,) = 2$. $T$ is not
Cohen-Macaulay exactly in the points of this finite set, therefore we must have
$V(a)\cap V(b) = \{\,m_1, n_1, m_2, n_2, m_3, n_3\,\}$. It follows that there
exist
linear forms $b_1$, $b_2$ and $b_3\in k{\lbrack}x_2,x_3,x_4 {\rbrack}$
such that $b = {b_1}{b_2}{b_3}$ and such that the 6-secant lines are given by
$L_i = \{\,x_0 = x_1 = b_i = 0\,\}$. For the last statement of the proposition
one computes $e(Z) =
e(Y) = -2$ and uses lemma 4.1, since neither $Z$ nor $Y$ are contained in a
quadric
hypersurface.$\slutt$\par

\proclaim Remark 4.8. {\rm The scheme $Z$ above is the union of the 5-secant
lines to $S$ and this
fits with the corresponding formula of Le Barz}.\par\bigskip

\proclaim Proposition 4.9. A smooth $K3$ surface $S$ of degree 10 in $\Pn 4$,
with $\pi=9$ and
one 6-secant can be bilinked $(4, 5)$ and $(4, 4)$ to a reducible surface $Z =
P\cup Q\cup T$,
where $T$ is a smooth cubic scroll, $Q$ is a smooth quadric surface cutting $T$
along two
rulings and $P$ is a plane cutting the scroll $T$ along its directrix, which in
turn is the
6-secant of the surface $S$. Moreover, $Z$ is the unique minimal element in its
even liaison class.
\par
$Proof.$ Let $V$ denote the unique quartic hypersurface containing $S\subset\Pn
4$,
cf. proposition $3.7$. Let also $H_{12}$ denote the hyperplane which contains
the
exceptional lines $E_1$ and $E_2$. By proposition 2.21, the residual
curve $H_{12}\cap S-E_1-E_2$ decomposes into two plane
quartic curves $A$ and $B$, with $A^2=B^2=0$ and $A\cdot B=2$. The planes
$\pi_1$ and $\pi_2$
spanned by $A$ and $B$, respectively, meet along the unique $6$-secant $L$ to
the surface $S$. Thus $\pi_i\subset V$, for $i=1,2$, and residual to them in
$V\cap H_{12}$
there exists a smooth quadric surface $Q$ meeting $S$ along $E_1\cup E_2$. By
proposition 3.7, the quintic
hypersurfaces containing $S$ cut out the surface outside the $6$-secant $L$, so
$S$ can be
linked in the complete intersection of $V$ and a general quintic hypersurface
$W$ to an irreducible
surface $Y$ of degree $10$, sectional genus $9$. From the liaison exact
sequences we get
$\cohi 0 Y 4=2$ this time, so we can link $Y$ on $V$ to a locally
Cohen-Macaulay
scheme $Z$, with $\deg\, Z=6$ and $\pi(Z)=1$. We identify in the sequel the
components of $Z$.
First we remark that, since $Q\cap Y=W\cap Q-E_1-E_2$ is a curve of type
$(3,5)$ on the quadric,
$Q$ is contained in all quartic hypersurfaces containing $Y$,
hence $Q$ is a component of $Z$. On another side, from lemma $1.4$ and
proposition $3.7$, we deduce that $\Delta(3)$ is a line, thus there exists a
plane $P$
such that all the hyperplanes $H\supset P$ have the
property that $\cohi 0 {H\cap S} 3 =1$. Any quartic containing $H\cap S$
depends
on the unique cubic through $H\cap S$, therefore $P$ lies inside $V$.\par

\proclaim Claim. $P$ meets the surface $S$ along a zero-dimensional scheme.\par

$Proof$ $of$ $the$ $claim.$ Assume that $P\cap S$ contains a curve $C$.
The general residual curve $H\cap S-C$, for $H\supset P$, is contained then in
an irreducible quadric surface inside $V$, since again any quartic surface in
$H$ containing $H\cap S$
depends on the unique cubic through $H\cap S$. Thus $V$ is singular along $P$,
and by linkage $P\cap Y$ is a curve of degree $5$.
It follows that $P$ is contained in all quartic
hypersurfaces containing $Y$, and thus it is a component of $Z$.
But the liaison exact sequences yield $Z\cap S=K=E_0\cup E_1\cup E_2$,
hence the intersection curve $C\subset P\cap S$ must then satisfy
$C\subset E_0\cup E_1\cup E_2$, which is a contradiction.$\slutt$\par\bigskip

Since $P$ meets $S$ only in points it follows that $P\cap Y$ contains a plane
quintic, and
thus as above, by Bezout, $P$ must be a component of the scheme $Z$.  The plane
$P$ contains
the $6$-secant line $L$ since
$\Delta(3)\subset\Delta(4)=L^{\scriptscriptstyle{\vee}}=\Pn 2$,
but it doesn't lie in $H_{12}={\rm span}_k(E_1,E_2)$. Therefore $P\cap Q=L\cap
Q=\lbrace{p_1, p_2}\rbrace$
and residual to $P\cup Q$ in $Z$ there is an irreducible surface $T$ of degree
$3$, which
contains the rational normal curve $E_0$. Moreover, since $Q\cap Y$ is a curve
of type
$(3,5)$ on the quadric, it follows that $T$ intersects $Q$ along two skew
lines, and hence
$T$ is a smooth cubic scroll meeting the quadric surface along two rulings. The
other claims of
the proposition follow from lemma $4.1$.$\slutt$\par\bigskip

{\bf 4.10. Construction of surfaces D.}
Let, as above, $T$ be a rational cubic scroll in $\Pn 4$ and let $Q$ be a
smooth quadric surface
cutting $T$ along two lines in its ruling, say $L_1$ and $L_2$ and consider
next a
plane $P$ passing through the directrix $L$ of $T$, cutting the scroll only
along this
line, and not contained in the hyperplane spanned by the quadric surface. Let
again
$Z=P\cup T\cup Q$.  $\deg Z=6$ and $\pi(Z)=1$ and $Z$ is locally Cohen-Macaulay
and a local
complete intersection except for the points $\lbrace p_i\rbrace=L\cap L_i$,
$i=\overline{1,2}$. We prove in the sequel that $Z$ can be backwards linked
($4,4$)
and ($4,5$) to a smooth $K3$ surface of type D. First a lemma.

\proclaim Lemma 4.11. \hfil\break
a) The scheme $X=T\cup Q$ is a degenerated elliptic quintic scroll in $\Pn 4$
and,
in particular, the homogeneous ideal $I_X$ is generated by $5$
cubics.\hfil\break
b) The homogeneous ideal $I_Z$ is generated by $10$ quartics and one quintic.
Moreover the quartics cut out the scheme $Z$ outside $L$ and the generic
quartic is singular along
and contains $L$ with multiplicity two.\par
$Proof.$ Let $H$ be the hyperplane spanned by the quadric surface $Q$.
The residual exact sequence
$$\esi T 2 X 3 {H\cap X,H} 3$$
remains exact after taking global sections since $\cohi 1 T 2=0$. Therefore
$\cohi 0 X 3=5$ and it suffices to check
whether $T$ is cut out by quadrics and $H\cap X$ by cubics. The former is clear
since
$T$ is defined by the minors of a $2\times 3$ matrix with linear entries, while
$\id {H\cap X,H} 3=\id {Q\cup L,H} 3\simeq\id {L,H} 1$ is clearly global
generated.
Moreover, taking cohomology in the above sequence we get
$$\cohi 1 X k=0\quad{\rm for\; all}\; k,\quad \cohi 2 X k=0\quad{\rm for }\quad
k\ne 0,
\quad{\rm and}\quad\hi 2({\cal I}_X)=\hi 2({\cal I}_{Q\cup L})=1,$$
thus Beilinson's spectral sequence gives a resolution of the form
$$\exactg {\Omega_{\Pn 4}^3(3)} {5\sstruct {\Pn 4}} {{\cal I}_X(3)}$$
and hence $X$ is a degenerated elliptic quintic scroll. Consider now the exact
sequence
$$\esi X {k-1} Z k {Z\cap H',H'} k$$
where $H'$ is a general hyperplane containing $P$. It remains also exact after
taking global
sections since $\cohi 1 X k=0$, for all $k$. Now $Z\cap H'=P\cup D\cup f_1\cup
f_2$,
where $D=Q\cap H'$ is a smooth conic and $f_1$, $f_2$ are rulings of $T$, thus
$\cohi 0 {Z\cap H'} k=\cohi 0 {D\cup f_1\cup f_2} {k-1}$ and we compute
$\cohi 0 {D\cup f_1\cup f_2} 2=0$,  $\cohi 0 {D\cup f_1\cup f_2} 3=\scohos 0
{\Pn 3} 3-
2\scohos 0 {\Pn 1} 3-\scohos 0 {\Pn 1} 6=5$ and $\scohid 0 {D\cup f_1\cup f_2}
4=16$.
Moreover, the homogeneous ideal $I_{D\cup f_1\cup f_2}$ is generated by the $5$
cubics and one
extra quartic and, since by $a)$ $I_X$ is generated by $5$ cubics, the first
assertion of $b)$
follows. For the second part it is enough to observe that the cubics in
$\sCohid 0
{D\cup f_1\cup f_2} 3$  vanishing on $L$ and  cut out, scheme-theoretically in
fact,
$D\cup f_1\cup f_2\cup L$.$\slutt$\par\bigskip

As a consequence of the above lemma  $Z$ can be linked in
the complete intersection of two quartic hypersurfaces to an irreducible
surface $Y$,
with $\deg\, Y=10$, $\pi(Y)=9$ which contains and is singular along $L$ and
which is smooth outside
this line. $Y$ can be further linked $(4,5)$ to a surface $S$ with the desired
invariants:
deg$S=10$, $\pi(S)=9$, and from the liaison exact sequences $p_g=1$, $q=0$. It
is
easily seen that $S$ is smooth outside $L$, for a general choice of the
linkages. In order to see the
behavior at the intersection with $L$ we'll work on the blow-up
$\widetilde{\Pn 4}$ of $\Pn 4$ along $L$, like in the proof of
$(4.5)$.\par\bigskip

Namely, keeping the same notations for the blow-up, if
$V_1$ and $V_2$ are two general quartic hypersurfaces containing $Z$, then, by
lemma $4.10$, they
have multiplicity two along $L$, thus $\overline{V_i}\equiv 2B$,
$i=\overline{1,2}$. On the
other hand the strict transforms of $P$ and $Q$ cut $E$ along
$\overline{P}\equiv
(B - F)(B - F)\equiv B^2 - 2BF$ and $\overline{Q}\equiv (B - F)2F\equiv 2BF$
respectively.
Also, as argued in $(4.5)$, the strict transform of $T$ cuts $E$ along
$\overline{T}\equiv B^2$.
It follows that, for a general $(4,4)$ linkage, the strict transform of $Y$ on
$\widetilde{\Pn 4}$
meets $E$ in a curve equivalent to $(2B)(2B) - (B^2 + (B^2 - 2BF) + 2BF)\equiv
2B^2$. A local
computation shows that the general quintic hypersurface containing $Y$ has
multiplicity one
along $L$. Therefore, for a general choice of the $(4,5)$ linkage, the strict
transform of $S$ on
$\widetilde{\Pn 4}$ will meet $E$ in a curve equivalent to $2B(B + 3F) -
2B^2\equiv 6BF$.
A Bertini argument shows now that for a general choice of the linkage, the
surface $S$ residual
to $Y$ is smooth. Moreover, since a curve of type $BF$ is blown down on $S$, it
follows that
$L$ is a $6$-secant line to $S$.\par\bigskip

To show that $S$ is indeed a $K3$ surface of type D we determine the one
dimensional components
in the intersections $S\cap Q$ and $S\cap T$. The liaison exact sequence  for
$Y$ gives
$P\cap Y\equiv 3H_P - K_P - L\equiv 5H_P$, $Q\cap Y\equiv 3H_Q - K_Q-\lbrack
Q\cap T\rbrack
\equiv 5l_1 + 3l_2$, where $l_1$ and $l_2$ denote the classes of the two
rulings of the
quadric, and $T\cap Y\equiv 3H_T - K_T - \lbrack T\cap(P\cup Q)\rbrack\equiv
3(C_0 + 2f) -
(-2C_0 - 3f) -C_0 - 2f\equiv 4C_0 + 7f$, with $C_0$ denoting the numerical
class of the
directrix $L$ on the scroll $T$ and $f$ the class of a ruling. The
one-dimensional components
of the intersection scheme $S\cap Z$ are the residuals (in term of conductor
ideals) of the
above curves in the complete intersection of $Z$ with the quintic hypersurface
used in the
linkage of $Y$ with $S$. Therefore, for a general choice of the linkage, $P$
cuts $S$ only in
points, $T$ cuts $S$ along a scheme whose one-dimensional part $K_1$ is
equivalent to
$5H_C -\lbrack T\cap Y\rbrack\equiv C_0 + 3f$ and $Q$ cuts $S$ along a curve
$K_2$ equivalent
to $5H_Q - \lbrack Q\cap T\rbrack\equiv 2l_2$ plus a zero dimensional scheme.
On the other
hand,  the scheme $K_1\cup K_2$ is exactly the canonical
divisor of $S$. Now the liaison exact sequence
$$\exactg {\id {Y\cup Z} 5} {\id {Y\cup P\cup T} 5} {\sstruct Q}$$
remains also exact after taking global sections, hence the quintics in $\sCohid
0 Y 5$
cut on $Q$ a linear system whose fixed part is exactly $Q\cap Y$. Therefore,
for a general
choice of the $(4,5)$ linkage, the curve $K_2\subset Z\cap S$ is reduced and
hence it is the
union of two skew lines, say $E_2$ and $E_3$, in the ruling of $Q$ containing
$L_1$ and
$L_2$. Eventually, the adjunction formula on $S$ yields $E_i^2 + KE_i = 2E_i^2
=
2p_a(E_i) - 2 =-2$, $i=\overline{2,3}$ and thus $E_2$ and $E_3$ are exceptional
lines on $S$.
It follows that $E_1$:=$K_1$ is a $(-1)$ quartic on $S$ and hence $S$ is, as
claimed,
a $K3$ surface of type D; i.e., embedded by a linear system of type
$$H=H_{min}-4E_1-E_2-E_3.$$

\par\bigskip

\proclaim Proposition 4.12. A smooth elliptic surface $S$ of degree 10 in $\Pn
4$, with $\pi=9$
can be bilinked $(4,5)$ and $(4,4)$ to a reducible surface $Z = M\cup Q$
where $M$ is a locally Cohen-Macaulay multiplicity four
structure on a plane $P$ given, after a suitable change of
coordinates and assuming $P = \{\,x_0 = x_1 = 0\,\}$, by the homogeneous ideal
$$I_M = {(\,{x}_0, {x}_1\,)}^3 + (\, g{x_0}^2 - f{x_0}{x_1},
h{x_0}^2-f{x_1}^2, h{x_0}{x_1} - g{x_1}^2\,)$$ with
$g, h, f \in k{\lbrack}x_2,x_3,x_4{\rbrack}$, homogeneous forms of degree 2
without common factors,  and where $Q$ is a smooth quadric surface in a
hyperplane non
containing $P$. Moreover, $Z$ is the unique minimal element in its even liaison
class.\par

$Proof.$ We keep the notations made in the proof of proposition 2.23 and the
lemmas thereafter. Namely, on $W$ the blow up, first of $\Pn 3$ at the vertex
$q$, then of the strict transform of $C_B$ and then of the strict transform of
the line $L\subset\Pn 6$, $h$ will denote the pullback of a plane, $E_L$ the
exceptional
divisor over $L$, $E_q$ the strict transform of the exceptional divisor
over $q$, and $E_B$ the strict transform of the exceptional divisor over $C_B$.
Then, for the strict transforms $\widetilde{S_1}$ and $\widetilde{H}$ of the
elliptic surface and of a hyperplane in $\Pn 4$ respectively, we have
$$\widetilde{S_1}\in |7h-2E_B-4E_q|$$
and
$$\widetilde{H}\in |4h-E_B-3E_q-E_L|.$$
By proposition 3.6, the quintic hypersurfaces containing $S$ cut out the
surface,
so $S$ can be linked in the complete intersection of the unique quartic $X$
containing it and a general quintic hypersurface to an irreducible surface
$Y$ of degree $10$ and sectional genus 9. Now, as in proposition 4.9,
$\cohi 0Y4=2$ so we can link further $Y$ in the complete intersection of $X$
with a general quartic hypersurface to a locally Cohen-Macaulay scheme $Z$
with $\deg Z=6$, $\pi(Z)=1$. We determine in the sequel the components of $Z$.
In terms of linear equivalence, for the strict transform $\widetilde{Z}$ of
$Z$, we have
$$\widetilde{Z}\equiv\widetilde{S_1}-\widetilde{H}\equiv 3h-E_B-E_q+E_L.$$
But a divisor in $|3h-E_B-E_q|$ is mapped down to $\Pn 3$ onto a cubic surface
containing the curve $C_B$, which thus necessarily coincides with the cone
$S_3$. It follows that in fact $\widetilde{Z}$ splits as
$$\widetilde{Z}\equiv (3h-E_B-3E_q)+2E_q+E_L,$$
and since $S_{3,W}\equiv 3h-E_B-3E_q$ gets contracted via the map to $\Pn 4$
onto the plane cubic $C$, while $E_q$ is mapped $2:1$ onto the plane $P$
spanned
by $C$, we deduce that $Z$ is the union of a multiplicity 4 structure on
the plane $P$ and of a quadric $Q$. The intersection between $S$ and $Q$ is
exactly the union of the 3 exceptional lines $E_i$, $i=\overline{1,3}$, so
in particular the quadric surface $Q$ is smooth. We identify in the sequel
the multiplicity four structure $M$ on the plane $P$. Let $Q^{\prime}$ be the
residual
quadric surface of $Q$ in the intersection of the quartic hypersurface $X$
with the hyperplane spanned by $Q$. It is easily seen that the union $Y\cup Q$
is linked in the complete intersection of $X$ and a sextic hypersurface to
the union $S\cup Q^{\prime}$. The scheme $S\cup Q^{\prime}$ is locally
Cohen-Macaulay since
any point on this scheme which is not of this type would lie necessarily
also on $Y$ and thus would be a base point of the whole linear system where
$Y$ is moving. But $|Y|$ has no base points on the surface $S$ since this
last one is cut out by quintic hypersurfaces on $X$. We deduce that $S\cup
Q^{\prime}$
is locally Cohen-Macaulay, and then by linkage $Y\cup Q$ is also of this type.
Finally, the above linkage between $Y$ and $Z$ implies that the multiplicity
four structure $M$ is a locally Cohen-Macaulay scheme. Moreover, the above
linkages yield $\cohi 0 M 2=\coh 0 {Y\cup Q}{K-H}=\cohi 0 {S\cup Q^{\prime}}
4-1=0$,
whence $M$ is a multiple structure of type $I_{X_9}$ in the list of [Ma],
i.e., of the type claimed in the statement of the proposition,  since
all the other multiplicity four Cohen-Macaulay structures in the list lie
on at least one hyperquadric.$\slutt$\par\bigskip

\proclaim Proposition 4.13. A smooth general type surface $S$ of degree 10 in
$\Pn 4$, with
$\pi=9$ can be bilinked $(4, 5)$ and $(4, 4)$ to a reducible surface $Z = T\cup
P_1\cup P_2\cup P_3$,
where $T$ is a smooth cubic Del Pezzo surface and $P_1$, $P_2$ and $P_3$ are
three planes
in general position, cutting $T$ along the three coplanar 6-secants of the
surface $S$. Moreover
$Z$ is a minimal element in the even liaison class of $S$.\par
$Proof.$ We study $\deh 3$. Recall from proposition $1.11$ that $\deh 4$ is the
union of
three planes meeting along a line $L$, corresponding to the pencil of
hyperplanes through
the $\Pn 2$ spanned by the $6$-secants. On the other side, $\deh 3\subset\deh
4$ is defined
by the minors of a $2\times 4$ matrix with linear entries, and moreover it
contains the line
$L$ by lemma $1.14$ and the remark before. Therefore, either $\deh 3$ contains
a plane, or
$\deh 3$ is a degenerated rational quartic curve, namely a union of lines
$L\cup L_1\cup L_2\cup L_3$ (some might coincide), such that no $3$ lines of
type $L_i$
are coplanar ( if they would be coplanar, the plane they span lies inside $\deh
3$).
Assume first that there exists a plane $\Pi$ such that $\Pi\subset\deh 3$, and
let $H$ be a hyperplane section of $S$ corresponding to a general point of
$\Pi$. Then
$\cohi 0H3=1$, $\cohi 1H3=3$ and $\cohi 0H4=4$, because otherwise $H\in\deh
{3,4}$ and
thus, by lemma 1.14 $H$, would be reducible, which is impossible. It follows
that the unique
quartic hypersurface $V$ containing $S$ (cf. $1.11$) splits off a plane in each
hyperplane
parameterized by $\Pi$. Hence $V$ is a cone with vertex the line $L$. But this
is impossible
since, by lemma 1.8, $S$ has no pencils of plane curves.
Therefore $\deh 3$ is a curve, and we may write $\deh 3=L\cup L_1\cup L_2\cup
L_3$, where
$L\cap L_i\ne\emptyset$, for all $i$. We observe also that two lines $L_i$ and
$L_j$ are
necessarily disjoint. Otherwise, the two corresponding dual planes
$L_i^{\scriptscriptstyle{\vee}}$
and $L_j^{\scriptscriptstyle{\vee}}$ would span only a hyperplane in $\Pn 4$,
and in that
hyperplane $V$ would decompose as $V=L_i^{\scriptscriptstyle{\vee}}\cup
L_j^{\scriptscriptstyle{\vee}}\cup Q$, where $Q$ is a quadric surface such that
$Q\supset{\rm span}(L_i,L_j)\cap S$, absurd. Therefore the planes $P_i:=
L_i^{\scriptscriptstyle{\vee}}$ meet pairwise only in points. As in proposition
4.9 one shows
that they meet $S$ only along zero-dimensional schemes, thus a Bezout argument
shows
that the $3$ planes $P_i$ are components of the bilinked scheme $Z$. Residual
to them,
there exist a scheme $T$ of degree $3$, which contains the intersection curve
$D=S\cap Z
\in |K_S|$. For a general choice of the linkage in the statement of the
proposition,
$D$ is an integral canonical curve of degree $6$ (cf. proof of $2.32$), and one
deduces easily
that $T$ is a smooth cubic Del Pezzo surface. The last statement in the
proposition follows from lemma $4.1$.$\slutt$.\par\bigskip

\proclaim Proposition 4.14. A smooth elliptic surface $S$ of degree 10 in $\Pn
4$, with $\pi=10$ can
be linked $(4,4)$ to a reducible surface $Z = T_1\cup T_2$, union of a smooth
Del Pezzo surface $T_1$
of degree 4 and a smooth quadric surface $T_2$ such that the
hyperplane of the quadric cuts $T_1$
along $G_1 + G_2 + F_1 + F_2$, consisting of
four lines with $G_1\cdot G_2 = F_1\cdot F_2 = 0$ and
such that $F_1$ and $F_2$ are members of one of the rulings
of $T_2$ and $G_1$ and $G_2$ meet
transversally the quadric. Furthermore the scheme $Z$ is minimal in its even
liaison class and can
be linked $(3,4)$ to a scheme $Y$ of degree 6, with $\pi=2$ which is minimal in
the even liaison
class of the surface $S$.\par
$Proof.$ We recall from proposition 2.25 that $S$ has only two exceptional
lines $E_1$ and $E_2$
and that $K\equiv F + E_1 + E_2$, with $|F|$ a base point free pencil defining
the elliptic
fibration. Consider the residual curve $D\equiv H - E_1 - E_2$ in the
intersection of $S$ with the
hyperplane spanned by the two exceptional lines. Then $D$ is a curve of degree
8 and arithmetic
genus 8, hence Riemann-Roch gives $\chi ({\cal O}_D (2)) = 9$. The cohomology
of the exact sequence
$$\es S {2H - D} S {2H} D {2H}$$ gives $h^1({\cal O}_D (2)) = 0$ since, by
lemma 1.4, $h^1({\cal O}_S (2)) = 0$
and also $h^2({\cal O}_S (2H - D)) = h^0({\cal O}_S (D - 2H + K)) = 0$.
Therefore taking global sections
in the exact sequence $$\eso D 2 {\Pn 3} 2 D 2 $$ shows that $D$ lies on a
quadric surface $T_2$. It
can't be a cone since then $D$ would have genus 9. Also $T_2$ is not reducible
since $S$ has no
6-secants and, by corollary 1.2, no plane curves of degree $\geq$ 5. It follows
that $T_2$ is a
smooth quadric in the hyperplane spanned by $E_1$ and $E_2$ and $D$ is a curve
of type $(3,5)$ on
it. Therefore, since $S$ is contained in a net of irreducible quartics, and by
Le Barz's formula
there are two 5-secants to $S$ meeting a general plane in $\Pn 4$, we may
characterize $T_2$ as the
union of the 5-secants to $S$. Any quartic hypersurface containing $S$ must
also contain $T_2$.\par
$S$ can be linked $(4,4)$ to a scheme $Z$ of degree 6, with $\pi=2$ which, by
the above discussion,
contains $T_2$ as a component. Therefore we get $Z = T_1\cup T_2$, with $T_1$ a
scheme of degree 4
such that $S\cap T_1 \in\mid 2H - F\mid$. On another side $T_1$ is the
projection from the improper node of the Del Pezzo surface $T\subset\Pn 5$ in
the proof of the proposition 2.34. So $T_1$ is smooth and $T_1\cap T_2$ is
of type $(2,0)$ on $T_2$. The configuration of lines in the proposition is
now clear, and the statements about linkage follow from lemma
4.1.$\slutt$\par\bigskip

{\bf 4.15. Construction of surfaces G. } We start as above. Let $T_1$ be a
smooth Del Pezzo
surface of degree $4$, and let $G_1+G_2+F_1+F_2$ be one of its hyperplane
sections, which
consists of four exceptional lines, such that say $G_1\cdot G_2=F_1\cdot
F_2=0$. Let now $T_2$ be
a smooth quadric surface in this hyperplane such that $F_1$ and $F_2$ are
members of one
its rulings, while $G_1$ and $G_2$ meet transversally $T_2$.\par

\proclaim Lemma 4.16. $T=T_1\cup T_2$ is a local complete intersection scheme,
it
lies on a pencil of reducible cubic hypersurfaces, and
its homogeneous ideal $I_T$ is generated by quartics.\par
$Proof.$ One uses the residual exact sequence
$$\esi {T_1} 3 T 4 {T\cap H} 4,$$
where $H$ is the hyperplane spanned by $T_2$ and argues as in lemma
$4.11$.$\slutt$
\par\bigskip

As a consequence of the previous lemma, we may link $T$ in the complete
intersection
of two general quartic hypersurfaces with a smooth surface $S$, with
invariants:
$\deg\, S=10$, $\pi=10$, and from the liaison exact sequences, $p_g=2$, $q=0$.
In particular, $S$ is by proposition $2.34$ a non-minimal proper elliptic
surface
with two exceptional lines. A closer look to the linkage, as in the
construction $4.10$,
shows in fact that the elliptic pencil of $S$ is given by the moving part of
the trace
on $S$ of the pencil of quadrics defining $T_2$.\par\medskip

\proclaim Remark 4.17. The scheme $T_2$ may be also characterized, in
accordance with
the corresponding formula of Le Barz, as the union of the $5$-secant lines to
$S$.\par\bigskip

\proclaim Proposition 4.18. A smooth general type surface $S$ of degree 10 in
$\Pn 4$, with $\pi=10$
can be linked $(4,4)$ to a reducible surface $Z = T_1\cup T_2$ where $T_1$ is a
cubic Del Pezzo
surface and $T_2$ is a degenerated cubic scroll, union of three planes having a
line in common and
such that $T_2$ cuts $T_1$ along a doubling of this line. Moreover the scheme
$Z$ is minimal in its
liaison class and can be linked $(3,4)$ to a scheme $Y$ of degree 6, with
$\pi=2$, which is minimal
in the even liaison class of $S$.\par
$Proof.$ It follows from proposition 2.34 that $S$ is embedded by $H\equiv 2K -
A_1 - A_2 - A_3$, with
$A_i$, for $i=\overline{1,3}$, disjoint $(-2)$-conics. Since $p_g = 3$ and
$K\cdot A_i = 0$, $i=\overline{1,3}$, by
taking global sections in the exact sequence $$\ess S {K - A_i - A_j} S K {A_i}
K {A_j} K$$ we get
that there exists a curve $D_k\in\mid K - A_i - A_j\mid$, for all triples
$\{\,i,j,k\,\} = \{\,1,2,3\,\}$.
It has degree 4 and genus 3, hence it is a plane curve. Now $D_1\cdot D_2 =
D_1\cdot D_3 = D_2\cdot D_3
= 2$ so the planes ${\Pi}_1$, ${\Pi}_2$ and ${\Pi}_3$ of $D_1$, $D_2$ and $D_3$
resp. meet pairwise
in lines. Since the three planes span all of $\Pn 4$ they must intersect along
a common line $L$. If
$L$ would lie on $S$ then, by corollary 1.2, it needs to be a component of all
curves $D_i$, $i=\overline{1,3}$,
whence the splitting $D_i = L + M_i$. But then $L\cdot M_i = 3$ and $L^2 +
L\cdot K = -2$ give
$K\cdot L - L\cdot A_j - L\cdot A_k = 1$, for all $\{\,i,j,k\,\} =
\{\,1,2,3\,\}$, and hence, summing up,
$6K\cdot L - 2(L\cdot A_1 + L\cdot A_2 + L\cdot A_3) = 3$, which is absurd.
Therefore $L$ doesn't lie
on $S$ so it must be the unique 6-secant for the surface. Furthermore, since
$S$ meets each of the
planes ${\Pi}_i$ in two points on $L$ outside the curve $D_i$,
$i=\overline{1,3}$, we get that each
plane  ${\Pi}_i$ has two pencils of 5-secant lines to $S$. Altogether we have
counted six 5-secants to
$S$ meeting a general plane in $\Pn 4$, and this fits with Le Barz's formula.
Therefore $T_2 = {\Pi}_1
\cup {\Pi}_2\cup {\Pi}_3$ is the union of the 5-secant lines to $S$.\par
By proposition 1.16, $S$ can be linked $(4,4)$ to a surface $Z$ of degree 6,
with $\pi=2$. $T_2$ is
contained in any quartic hypersurface containing $S$ so the planes ${\Pi}_1$,
${\Pi}_2$ and ${\Pi}_3$
will necessarily be components of $Z$. Let $T_1$ be the residual component:
$$Z = T_1\cup {\Pi}_1\cup {\Pi}_2\cup {\Pi}_3.$$ Thus $T_1$ has degree 3. The
liaison exact sequence
$$\es S K S 3 {S\cap Z} 3$$ gives $S\cap Z\equiv 3H - K\equiv H + D_1 + D_2 +
D_3\equiv H + \sum_{i=1}^3
{\Pi}_i\cap S$, hence $T_1\cap S\equiv H$ on $S$ and $T_1$ is contained in a
hyperplane. Now
$T_1\cap (S\cup {\Pi}_1\cup {\Pi}_2\cup {\Pi}_3)\equiv 3H_{T_1} - K_{T_1}\equiv
4H_{T_1}$, hence $T_1
\cap S$ is linked in the intersection of a quartic and the cubic $T_1$ to a
curve $F$ of degree 2 and
arithmetic genus $p_a(F) = -2$. This means that $F$ is a double structure on a
line $L^{\prime}$ on the
cubic $T_1$. In fact the line $L^{\prime}$ must coincide with the line $L$
since $T_2$ has no components
in the hyperplane of $T_1$ and the doubling of $L^{\prime}$ is contained in
${\Pi}_1\cup {\Pi}_2\cup {\Pi}_3$. It is also easily seen that, for a general
choice
of the linkage, $T_1$ is a smooth cubic Del Pezzo surface meeting
$T_2$ as stated in the proposition.
The rest follows from lemma $4.1$.$\slutt$\par\bigskip

\proclaim Remark 4.19. Using, for instance, the above description of $Z$ one
can show
that, for generic choices in the second linkage, $Y$ is also a reducible
surface
$F_1\cup F_2$, where $F_1$ is a smooth symmetric Castelnuovo surface and $F_2$
is a general plane
cutting $F_1$ along the singular line of the rank 3 hyperquadric containing the
Castelnuovo surface.
This line is exactly the unique 6-secant to the surface $S$.\par\bigskip

{\bf 4.20. Construction of surfaces H.}  As for surfaces of type $G$,
we start with a Del Pezzo surface of degree 3 and three planes
forming a degenerate cubic scroll, such that they all intersect the
hyperplane of the Del Pezzo along the same line $L$ on the Del Pezzo. The
union of these four surfaces is easily seen to be cut out by quartics, and
furthermore performing linkage on the blow up of $\Pn 4$ along the line
$L$, like in proposition $4.5$, one shows that
the union of the four surfaces is linked $(4,4)$ to a smooth
surface of degree 10, for which $L$ becomes a 6-secant line.\par

\vfil\eject

{\hbox {}}\bigskip
\centerline{\third   References}
\bigskip\bigskip
{\leftskip=6pt \rightskip=-6pt

\item{[Al]} Alexander, J.: Surfaces rationelles non-speciales dans $\Pn 4 $.
Math. Z. {\bf 200}, 87-110 (1988).

\item{} Alexander, J.: Speciality one  rational surfaces in
$\Pn 4$. LNS {\bf 179}, 1-23, London Math. Soc., Cambridge Univ. Press (1992).

\item{[ADHPR]} Aure, A., Decker, W., Hulek, K., Popescu, S., Ranestad, K.:
The Geometry of Bielliptic Surfaces in $\Pn 4$, Int. J. Math. (to appear).

\item{[AR]} Aure, A., Ranestad, K.: The Smooth Surfaces of Degree $9$ in
$\Pn 4$. LNS {\bf 179}, 32-46, London Math. Soc., Cambridge Univ. Press (1992).

\item{[Ba]} Baker, H. F., Principles of Geometry {\bf VI}, Cambridge Univ.
Press, p.271 ff. (1910).

\item{[BHM]} Barth, W., Hulek, K., Moore, R. Degenerations of Horrocks-Mumford
surfaces. Math. Ann. {\bf 277}, 735-755 (1987).

\item{[BaM]} Barth, W., Moore, R.: Geometry in the space of Horrocks-Mumford
 surfaces. Topology {\bf 28}, 231-245 (1989).

\item{[Mac]} Bayer, D., Stillman, M.:  Macaulay: A system for computation in
algebraic geometry and commutative algebra. Source and object code available
for Unix and Macintosh computers. Contact the authors, or download from
{\bf  zariski.harvard.edu} via anonymous ftp.

\item{[Bei]} Beilinson, A.: Coherent sheaves on $\Pn N$ and problems of linear
algebra. Funct. Anal. Appl. {\bf 12}, 214-216 (1978).

\item{[BoM]} Bolondi, G., Migliore, J.: The structure of an even liaison class.
Trans.
AMS. {\bf 316}, No. {\bf 1}, 1-38 (1989).

\item{[BF]} Braun, R., Fl\o ystad, G.: A bound for the degree of surfaces in
$\Pn 4$ not of general type, Preprint Bergen, Bayreuth (1993).

\item{[Br]} Brivio, S.: Smooth Enriques surfaces in $\Pn 4$ and exceptional
bundles, Math. Z. {\bf 213}, 509-521 (1993).

\item{[Co]} Comessati, A., Sulle superficie di Jacobi, Tipografia della
R. Accademia dei Lincei, Rome (1919).

\item{[D]} Decker, W.: Stable rank $2$ vector bundles with Chern-classes
$c_1=-1,
c_2=4$. Math. Ann. {\bf 275}, 481-500 (1986).

\item{[DES]} Decker, W.,Ein, L., Schreyer, F.-O.: Construction of surfaces in
$\Pn 4$. J. of Algebraic Geometry {\bf 2}, 185-237 (1993).

\item{[EP]} Ellia, Ph., Peskine, Ch.: Groupes de points de $\Pn 2$; caract\`ere
et
position uniforme. Algebraic Geometry. Proceedings, L'Aquila 1988, LNM. {\bf
1417},
111-116 (1990).

\item{[ElP]} Ellingsrud, G., Peskine, Ch.: Sur les surfaces lisses de $\Pn 4$.
Inv. Math. {\bf 95}, 1-11 (1989).

\item{[Ha]} Hartshorne, R.: Algebraic geometry. Springer (1977).

\item{[HM]} Horrocks, G., Mumford, D.: A rank 2 vector bundle on $\Pn 4$ with
$15000$ symmetries. To\-po\-lo\-gy {\bf 12}, 63-81 (1973).

\item{[H]} Hulek, K.: Geometry of the Horrocks-Mumford bundle. Algebraic
Geometry,
 Bowdoin 1985. Proc. Symp. Pure Math.
{\bf 46}, (2), 69-85 (1987).

\item{[HKW]} Hulek, K., Kahn, C., Weintraub, S.: Moduli Spaces of Abelian
Surfaces: Compactification, Degenerations, and
       Theta Functions.  Walter de Gruyter, Berlin (1993).

\item{[HL]} Hulek, K., Lange, H.: Examples of abelian surfaces in $\Pn 4$.
J. Reine Angew. Math.  {\bf 363}, 201-216 (1985).

\item{[Io]} Ionescu P.: Embedded projective varieties of small invariants I.
Proc.
of the week of Algebraic Geometry, Bucharest (1982), LNM. {\bf 1056}, 142-186
(1984).

\item{} Ionescu P.: Embedded projective varieties of small invariants II.  Rev.
Roum. Math. Pures Appl. {\bf 31}, 539-544 (1986).

\item{} Ionescu, P.:  Embedded projective varieties of small invariants
III. Proceedings of the l'Aquila conference. LNM. {\bf 1417}, 138-154 (1990).

\item{[K]} Kleiman, S.: Geometry on grassmanians and applications to splitting
bundles and smoothing cycles. IHES, Publ. Math. {\bf 36}, 282-298 (1969).

\item{[L]} Lange, H.: Embeddings of Jacobian surfaces in $\Pn 4$.
J. Reine Angew. Math.  {\bf 372}, 71-86 (1986).

\item{[LR]} Lazarsfeld, R, Rao, P.: Linkage of general curves of large degree.
Algebraic
Geometry - open problems, Ravello 1982, LNM. {\bf 997}, 267-289, Springer 1983.

\item{[LB]} Le Barz, P.: Formules pour les multisecantes des surfaces. C.R.
Acad.
Sc. Paris, t.{\bf 292} Serie I, 797-799 (1981).

\item{[Ma]} Manolache, N.: Codimension Two linear Varieties with Nilpotent
Structures. Math. Z. {\bf 210}, 573-580 (1992).

\item{[MDP]} Martin-Deschamps, M., Perrin, D.: Sur la classification des
courbes
gauches. Ast\'erisque {\bf 184-185} (1990).

\item{[Mo]} Moishezon, B.: Complex Surfaces and Connected Sums of Complex
Projective Planes. LNM. {\bf 603} (1977).

\item{[Ok]} Okonek, Ch.: Moduli reflexiver Garben und Fl\"achen von kleinem
Grad
in $\Pn 4$. Math. Z. {\bf 184}, 549-572 (1983).

\item{} Okonek, Ch.: \"Uber $2$-codimensionale Untermannigfaltigkeiten vom
Grad $7$ in $\Pn 4$ und $\Pn 5$. Math. Z. {\bf 187}, 209-219 (1984).

\item{} Okonek, Ch.: Fl\"achen vom Grad 8 in $\Pn 4$. Math. Z. {\bf 191},
207-223 (1986).

\item{[PS]} Peskine, Ch., Szpiro, L.: Liaison des vari\'et\'es alg\'ebriques I.
Invent.
Math. {\bf 26}, 271-302 (1974).

\item{[Po]} Popescu, S.: On Smooth Surfaces of Degree $\geq$ 11 in the
Projective
Fourspace. Thesis, Univ. des Saarlandes (1993).

\item{[R]} Ramanan, S.: Ample divisors on abelian surfaces. Proc. London
Math. Soc. {\bf 51}, 231-245 (1985).

\item{[Ra]} Ranestad, K.:On smooth surfaces of degree ten in the projective
fourspace. Thesis, Univ. of Oslo  (1988).

\item{} Ranestad, K.:Surfaces of degree 10 in the projective
fourspace. Symp. Math. Academic Press Vol.{\bf XXXII}, 271-307 (1991).

\item{[Ro]} Roth, L.: On the projective classification of surfaces. Proc.
London
Math. Soc. {\bf 42}, 142-170 (1937).

\item{[S]} Serrano, F.: Divisors on bielliptic surfaces and embeddings in $\Pn
4$.
 Math. Z. {\bf203}, 527-533 (1990).

\item{[Se]} Severi, F.: Intorno ai punti doppi impropri di una superficie
generale dello spazio a quattro dimensioni, e a suo punti tripli apparenti.
Rend. Circ. Mat. Palermo {\bf 15}, 33-51 (1901).
\par}

\bye